%% file: main.tex
\documentclass[aps,pre,groupedaddress,superscriptaddress,showpacs,nofootinbib,notitlepage]{revtex4-2}

\usepackage{makecell} % for splitting cell for eta_3
\usepackage{xurl}
\usepackage[hidelinks]{hyperref}
\usepackage{booktabs}       % professional-quality tables
\usepackage{amsfonts}       % blackboard math symbols
\usepackage{amssymb}
\usepackage{mathtools}
\usepackage{array}
\usepackage{ragged2e}
\usepackage{threeparttable}
\usepackage{nicefrac}       % compact symbols for 1/2, etc.
\usepackage{microtype}      % microtypography
\usepackage{amsmath}
\usepackage{xcolor,colortbl}
\usepackage{hhline}
\usepackage{bbm}
\usepackage[draft]{changes}
\usepackage{tabularx}
\usepackage{amsthm}
\usepackage{changes}
\usepackage{comment}
\usepackage{pifont}
\usepackage{color,soul}
\usepackage{xspace}
\usepackage{tikz}
\usetikzlibrary{positioning, shapes.geometric, arrows.meta}
\usepackage{multirow}
\usepackage{etoolbox}
\usepackage{rotating}
\usepackage{tcolorbox}
\patchcmd{\section}
  {\centering}
  {\raggedright}
  {}
  {}
\patchcmd{\subsection}
  {\centering}
  {\raggedright}
  {}
  {}

\newcolumntype{P}[1]{>{\centering\arraybackslash}p{#1}}
\newcolumntype{C}[1]{>{\centering\arraybackslash}p{#1}}
\newcolumntype{L}[1]{>{\raggedright\arraybackslash}p{#1}}

\def\ie{\textit{i.e.}}
\def\eg{\textit{e.g.}}

\makeatletter
\newcommand{\lowersim}[2]{%
  \sbox\z@{$#1<$}%
  \raisebox{-\dimexpr\height-\ht\z@}{$\m@th#1#2$}%
}
\makeatother

\makeatletter
\DeclareRobustCommand\onedot{\futurelet\@let@token\@onedot}
\def\@onedot{\ifx\@let@token.\else.\null\fi\xspace}

\definecolor{tabblue}{RGB}{31, 119, 180}
\definecolor{taborange}{RGB}{255, 127, 14}
\definecolor{tabgreen}{RGB}{44, 160, 44}
\definecolor{tabred}{RGB}{214, 39, 40}
\definecolor{tabpurple}{RGB}{148, 103, 189}
\definecolor{tabbrown}{RGB}{140, 86, 75}
\definecolor{tabpink}{RGB}{227, 119, 194}
\definecolor{tabgray}{RGB}{127, 127, 127}
\definecolor{tabolive}{RGB}{188, 189, 34}
\definecolor{tabcyan}{RGB}{23, 190, 207}

\definecolor{grey}{gray}{0.5}

\begin{document}
\title{Learning dynamical systems with biochemically informed neural ordinary differential equations}

\author{Luis L. Fonseca}
\thanks{\href{mailto:llfonseca@medicine.ufl.edu}{llfonseca@medicine.ufl.edu}}
\affiliation{Laboratory for Systems Medicine, Department of Medicine, University of Florida, Gainesville, FL, USA}
\author{Reinhard C. Laubenbacher}
\thanks{\href{mailto:reinhard.laubenbacher@medicine.ufl.edu}{reinhard.laubenbacher@medicine.ufl.edu}}
\affiliation{Laboratory for Systems Medicine, Department of Medicine, University of Florida, Gainesville, FL, USA}
\author{Lucas B\"{o}ttcher}
\thanks{\href{mailto:l.boettcher@fs.de}{l.boettcher@fs.de}}
\affiliation{Department of Computational Science and Philosophy, Frankfurt School of Finance and Management, Frankfurt am Main, 60322, Germany}
\affiliation{Laboratory for Systems Medicine, Department of Medicine, University of Florida, Gainesville, FL, USA}
\date{\today}
\begin{abstract}
Ordinary differential equation models of biochemical reactions are often formulated as stoichiometric systems in which the dynamics arise from a collection of interacting processes. A central challenge is that the functional form of each process is rarely known a priori and may be difficult to infer from data. We propose biochemically informed neural ordinary differential equations (BINODEs), a neural-ODE framework that retains the stoichiometric structure of mechanistic models while representing individual processes by neural networks. In BINODEs, the outputs of neural network processes (NNPs) are mapped to state derivatives through a linear layer analogous to a stoichiometric matrix. This architecture allows biological side information, such as process-specific inputs, sign constraints, and monotonicity assumptions, to be built directly into the model. We characterize the approximation properties of $\mathrm{NNP}$s for several standard biochemical rate laws and show that the proposed framework recovers both trajectories and process-level structure in Monod, Lotka--Volterra, pharmacokinetic, and ultradian endocrine models. These results suggest that BINODEs offer a useful compromise between mechanistic interpretability and data-driven flexibility for modeling partially known biochemical or biological dynamical systems.
\end{abstract}
\maketitle
\section{Introduction}
In biochemistry and related fields, mechanistic ordinary differential equation (ODE) models are widely used to describe biochemical systems. In such models, each differential equation describes the time evolution of one state variable (\eg, the concentration of a species). Mechanistically, these changes are modeled as the sum of the contributions from all processes affecting that variable. That is,
\begin{equation}
    \frac{\mathrm{d}X}{\mathrm{d}t}=N\, V(X)\,,    
    \label{eq:NP}
\end{equation}
where $X \equiv X(t) = (x_1(t), \dots, x_n(t))^\top \in \mathbb{R}^n$ is the system state at time $t$, with $x_i(t)$ denoting state variable $i\in\{1,\dots,n\}$. The matrix $N \in \mathbb{R}^{n \times k}$ is the stoichiometric matrix, whose entries $n_{ij}$ describe the effect of process $j \in \{1,\dots,k\}$ on $x_i$. The vector $V(X) = (v_1(X), \dots, v_k(X))^\top \in \mathbb{R}^k$ consists of (usually nonnegative) process rates $v_j(X)$.

Common mechanistic process representations include the Henri--Michaelis--Menten~\cite{michaelis1913kinetik,henri1903lois,henri2006theorie,cornish2015one} and mass-action~\cite{waage1864studier,voit2015150} laws. Although the mass-action law was originally formulated by Waage and Guldberg~\cite{waage1864studier} for chemical equilibria in solution, it was later found to effectively describe homogeneous populations in ecology~\cite{lotka1925elements,volterra1927variazioni} and epidemiology~\cite{kermack1927contribution}. More recently, both mass-action and frequency-based transmission models, in which an infected individual is assumed to have a fixed number of contacts over a given time period regardless of population density, have been studied in the context of pair transmission models~\cite{wylie2021uniformly}.
\begin{sidewaystable*}[p]
\centering
\small
\renewcommand{\arraystretch}{1.25}
\setlength{\tabcolsep}{5pt}

\begin{threeparttable}
\caption{Commonly used process representations in biochemistry and biology.}
\label{tab:processes}

\begin{tabularx}{\textheight}{@{}
    >{\RaggedRight\arraybackslash}p{3.8cm}
    >{\centering\arraybackslash}m{8.8cm}
    >{\centering\arraybackslash}p{2.0cm}
    >{\RaggedRight\arraybackslash}X
@{}}
\toprule
\textbf{Process type} & \textbf{Mathematical form} & \textbf{Symmetric} & \textbf{Additional comments} \\
\midrule

Michaelis--Menten
& $\displaystyle \frac{V_{\max}S}{S + K_{\mathrm m}}$
& n/a
& $V_{\max}$ is the maximum reaction rate, and $K_{\mathrm m}$ is the Michaelis--Menten constant for substrate $S$, \ie, the substrate concentration for which the reaction rate is half-maximal \cite{michaelis1913kinetik,henri1903lois,henri2006theorie,cornish2015one}. \\

\makecell[l]{Random Bi--Bi\\$A+B \xrightarrow{E} P+Q$}
&
$
\displaystyle
\frac{
  V^+_{\max}\,\dfrac{A}{K_{\rm i}^A}\,\dfrac{B}{K_{\rm m}^B}
}{
  \begin{aligned}
  &1+\dfrac{A}{K_{\rm i}^A}+\dfrac{B}{K_{\rm i}^B}+\dfrac{P}{K_{\rm i}^P}+\dfrac{Q}{K_{\rm i}^Q}\\
  &\quad+\dfrac{A}{K_{\rm i}^A}\dfrac{B}{K_{\rm m}^B}
  +\dfrac{P}{K_{\rm m}^P}\dfrac{Q}{K_{\rm i}^Q}
  \end{aligned}
}
$
& no
& $V^+_{\max}$ is the maximum rate of the forward reaction, and $K_{\rm i}^X$ and $K_{\rm m}^X$ are the Michaelis--Menten constants associated with each substrate or product $X$.\\
\addlinespace[12pt]

Power-law
& $\displaystyle \alpha \prod_i X_i^{g_i}$
& yes
& $\alpha$ is the rate constant, representing the proportionality between the reaction rate and the product of concentrations $\prod_i X_i^{g_i}$; $g_i$ are the kinetic orders of each $X_i$~\cite{savageau1969iochemical,savageau1970biochemical,voit2013biochemical}. \\

Lin--log
& $\displaystyle \frac{v}{J^0}=\frac{e}{e^0}\left(1+\sum_i \epsilon_i \ln\frac{x_i}{x_i^0}\right)$
& yes\tnote{a}
& $v$ is the reaction rate, $e$ the enzyme level, $x_i$ the level of metabolite $i$, $\epsilon_i$ the elasticities associated with metabolite $i$, and $J^0$, $e^0$, and $x_i^0$ the reaction flux, enzyme level, and metabolite levels at reference state $0$, respectively~\cite{visser2003dynamic}. \\

Convenience kinetics
&
$
\displaystyle
\begin{aligned}
&\frac{
  V^+_{\max}\prod_i \widetilde{S}_i^{\alpha_i}
}{
  \prod_i\left(\sum_{m=0}^{\alpha_i}(\widetilde{S}_i)^m\right)
  +\prod_j\left(\sum_{m=0}^{\beta_j}(\widetilde{P}_j)^m\right)
}\\[1mm]
&\qquad\times
\prod_l\frac{A_l}{k^{A_l}+A_l}
\prod_n\frac{k^{I_n}}{k^{I_n}+I_n}
\end{aligned}
$
& yes\tnote{b}
& $\widetilde{X}=X/K_{\rm m}^X$, $K_{\rm m}^X$ is the Michaelis--Menten constant for $X$, $S_i$ is the $i$th substrate, $P_j$ the $j$th product, $A_l$ the $l$th activator, $I_n$ the $n$th inhibitor, $k^{A_l}$ is the activation constant for the $l$th activator, and $k^{I_n}$ is the inhibition constant for the $n$th inhibitor. \\
\addlinespace[12pt]

Gene transcription
& $\displaystyle k_{\rm t} n \frac{1+\sum_j k_{\mathrm{a},j}C_j}{1+\sum_l k_{\mathrm{r},l}C_l}$
& yes\tnote{c}
& $k_{\rm t}$ is the basal transcription rate, $n$ is the copy number of the gene of interest, $k_{\mathrm{a},j}$ is the strength of activation by gene $j$, $C_j$ the concentration of gene $j$, $k_{\mathrm{r},l}$ the strength of inhibition by gene $l$, and $C_l$ the concentration of gene $l$~\cite{acon2021myc}. \\

\bottomrule
\end{tabularx}

\begin{tablenotes}[flushleft]
\footnotesize
\item[a] The lin--log representation is symmetric with respect to metabolites $x_i$, but not with respect to the enzyme level.
\item[b] All substrates, products, activators, and inhibitors are symmetric within each category, but not across categories.
\item[c] All activators and inhibitors are symmetric within each category, but not across categories.
\end{tablenotes}

\end{threeparttable}
\end{sidewaystable*}

The formalism developed by Henri, Michaelis, and Menten for single-substrate enzyme-catalyzed reactions was later generalized to include multiple substrates and regulatory factors~\cite{cleland1975partition,king1956schematic,briggs1925note}. These generalizations also lead to the Hill equation when applied to proteins with multiple binding sites~\cite{hill1910possible,gesztelyi2012hill}. Yet these representations present a recurring trade-off: they can be difficult to infer from data, they rely on assumptions about the catalytic mechanism, and they often yield non-symmetric functions\footnote{We call a function $V$ \emph{symmetric} if it is invariant under simultaneous permutations of the state variables and their associated parameters. That is, for any permutation $\pi$ of $\{1,\dots,n\}$, $V(x_1,\dots,x_n \mid \theta_1,\dots,\theta_n)
=
V(x_{\pi(1)},\dots,x_{\pi(n)} \mid \theta_{\pi(1)},\dots,\theta_{\pi(n)})$, where $\theta_i$ denotes the collection of parameters associated with $x_i$.}, which can complicate model analysis. To address these limitations, more phenomenological process representations have been developed. Prominent examples include the power-law formalism, introduced in biochemical systems theory~\cite{savageau1969iochemical,savageau1970biochemical,voit2013biochemical}, and the lin-log formalism~\cite{visser2003dynamic}, developed for metabolic control analysis~\cite{kacser1973control,heinrich1974linear}. Other process representations proposed in biochemistry include the convenience rate law~\cite{liebermeister2006bringing}, the modular rate law~\cite{liebermeister2010modular}, and the saturable and cooperative formalism~\cite{sorribas2007cooperativity}. 

Similar variation in modeling approaches arose in other biological fields. In microbiology, differences in substrate uptake kinetics have led to a range of empirical growth models (see \cite{muloiwa2020comparison}). When microbial growth is assumed to be proportional to substrate uptake, the resulting growth law reflects the uptake kinetics. For example, Michaelis--Menten uptake leads to Monod-type growth kinetics~\cite{monod1958recherches, monod1950technique}, Haldane uptake leads to the Haldane growth model~\cite{haldane1965enzymes,andrews1968mathematical}, and Hill-type uptake leads to the Moser model~\cite{moser1958dynamics}. Likewise, in ecology, the Holling type II equation~\cite{cs1959components, holling1959some} is structurally similar to the Michaelis--Menten and Monod equations, even though it is derived in a mechanistically different way~\cite{dawes2013derivation}.

Taken together, these examples illustrate both the breadth of available process representations and the absence of a unified, flexible framework. Different process models have been introduced to describe catalysis, transport, regulation, predation, and many other biological processes, often balancing mechanistic fidelity, interpretability, tractability, and empirical adequacy in different ways. Table~\ref{tab:processes} summarizes common process representations used in computational biochemistry, systems biology, ecology, and related disciplines, and highlights key structural features (\eg, symmetry, saturability, and monotonicity). In Table~\ref{tab:processes_appendix} of Appendix~\ref{app:common_processes}, we summarize further commonly used process representations in biochemistry and biology. This variety, together with the limitations of individual rate laws, motivates a flexible representation that can incorporate biochemical structure without committing to a single closed-form expression. Rather than choosing among the many mathematical expressions used to describe biochemical processes, we propose biochemically informed neural ordinary differential equations (BINODEs) and suggest using neural network processes (NNPs) for the representation of biological processes. Our approach allows us to incorporate different types of prior information, including knowledge about the inputs and outputs of specific processes and monotonicity with respect to selected input variables.

The approach pursued in this paper is related to several recent developments at the interface of machine learning and mechanistic modeling in biochemistry and related fields. For example, to alleviate the parametric restrictions imposed by conventional ODE models, neural ODEs incorporating Hill--Langmuir kinetics have been proposed for modeling gene regulatory networks~\cite{hossain2024biologically}. In~\cite{philipps2024universal}, neural ODE functions have been used to capture unknown effects in the evolution of biological systems, including glycolysis and STAT5 dimerization. Similar hybrid approaches have been proposed in \cite{de2025physiology}, \cite{grigorian2024hybrid}, and \cite{thöni2025modeling}, and have been applied to the human cardiovascular system, Michaelis--Menten kinetics, and chemical reaction networks, respectively.

All of these hybrid mechanistic-neural approaches share the common feature of integrating existing knowledge into the formulation of artificial neural networks, thereby introducing an inductive bias that can facilitate training and improve generalization and consistency with known physical constraints (\eg, in control tasks~\cite{bottcher2026control}). This class of methods is known as dynamics-informed or physics-informed learning and, more generally, as learning with side information~\cite{DBLP:conf/l4dc/AhmadiK20,ahmadi2023learning}. In contrast to prior approaches that encode specific kinetic laws, we focus on incorporating more general biochemical structure, such as process-level input-output relations and monotonicity constraints.

The proposed BINODEs have potential applications as surrogate models for complex biomedical systems involving stochastic effects, spatial heterogeneity, and multiscale dynamics. By integrating biochemical structure directly into the neural representation of processes, this approach complements existing surrogate models, including traditional surrogates~\cite{fonseca2025optimal}, neural ODE surrogates~\cite{fronk2023interpretable,bottcher2025control}, and neural stochastic differential equation (SDE) surrogates~\cite{dietrich2023learning,zhang2025reconstructing}.

BINODEs provide a useful middle ground between conventional mechanistic ODE models and fully black-box neural ODEs. Compared with traditional rate-law models, they avoid committing to a potentially misspecified closed-form process representation; compared with generic neural ODEs, they retain an explicit stoichiometric decomposition and allow biochemical side information to be imposed directly. 

The remainder of the paper is organized as follows. We first introduce neural network processes ($\mathrm{NNP}$s) as modular process models, then show how NNPs compose into BINODEs, and finally apply BINODEs to bioreactor, Lotka--Volterra, pharmacokinetics, and ultradian endocrine dynamics. In the pharmacokinetics and ultradian endocrine examples, we compare our approach with that developed in~\cite{AI-Aristotle} and show that relatively small BINODEs can recover underlying state-dependent process representations, yielding an autonomous and more mechanistically interpretable description of the dynamics than time-dependent parameterizations. We conclude with a discussion of our results and outline promising directions for future work.
\section{Neural network process representation}
Stoichiometric systems of the form given in Eq.~\eqref{eq:NP} provide a common framework for ODE-based models in the biological sciences. Their structure is determined by the stoichiometric matrix, while each process rate is specified by a separate functional representation derived from first principles (\eg, mass-action or Michaelis--Menten kinetics) or chosen for mathematical convenience, structural properties, or empirical adequacy. In this section, we introduce neural network processes ($\mathrm{NNP}$s) as general approximators of process rates and examine their approximation properties.
\subsection{Formulation}
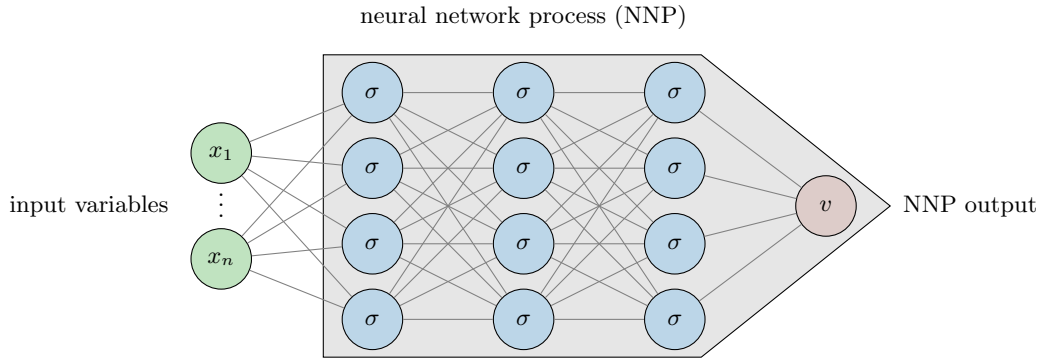
\begin{figure}
    \centering
    \input{nnp.tikz}
    \caption{Schematic representation of a neural network process (NNP), implemented as a feedforward network with input $X=(x_1,\dots,x_n)^\top$ and scalar output $v(X)$. We use $\sigma(\cdot)$ to denote an activation function, such as the rectified linear unit ($\rm ReLU$) or exponential linear unit ($\rm ELU$).}
    \label{fig:NNP}
\end{figure}
Each $\mathrm{NNP}$ maps the relevant system variables $x_1,\dots,x_n$ to a single scalar output representing the corresponding process rate, using a fully connected feedforward network (see Fig.~\ref{fig:NNP}). In the hidden layers, we primarily use the exponential linear unit ($\rm ELU$) activation function, which helps mitigate the ``dead ReLU problem'' (\ie, zero outputs and vanishing gradients for negative pre-activations) often associated with the rectified linear unit (ReLU)~\cite{DBLP:conf/acssc/DouglasY18}. For the output node, we use a ReLU or softplus activation to enforce nonnegativity of the process rates when required.

Representing biological processes with $\mathrm{NNP}$s offers two main advantages. First, under standard conditions, feedforward neural networks are universal function approximators~\cite{hornik1989multilayer,DBLP:conf/nips/LuPWH017,DBLP:conf/colt/Telgarsky16}. Second, multiple $\mathrm{NNP}$s can be combined into a neural ODE that preserves the stoichiometric structure of Eq.~\eqref{eq:NP}.

In analogy with the power-law representations used in biochemical systems theory~\cite{savageau1969iochemical,savageau1969biochemical,savageau1970biochemical}, we define the process rate associated with an $\mathrm{NNP}$ as
\begin{equation} 
    v(X)=\mathrm{NNP}(X;\theta)\,,    
    \label{eq:NNP}
\end{equation}
where $\theta\in\mathbb{R}^m$ denotes the $\mathrm{NNP}$ parameters.

The corresponding power-law process representation is
\begin{equation} 
    v(X)=\alpha \prod_i x^{g_i}_i\,,  
    \label{eq:power_law_process}
\end{equation}
where $\alpha$ and $g_i$ denote the rate constant and kinetic orders, respectively.

As in power-law formulations, which may depend on all or only a subset of the system variables (with $g_i=0$ for variables that do not influence the process), an $\mathrm{NNP}$ can be constructed as a function of all variables or only the relevant subset.
\subsection{Approximation properties}
\label{sec:characterizing_nnps}
\begin{figure}
    \centering
    \includegraphics{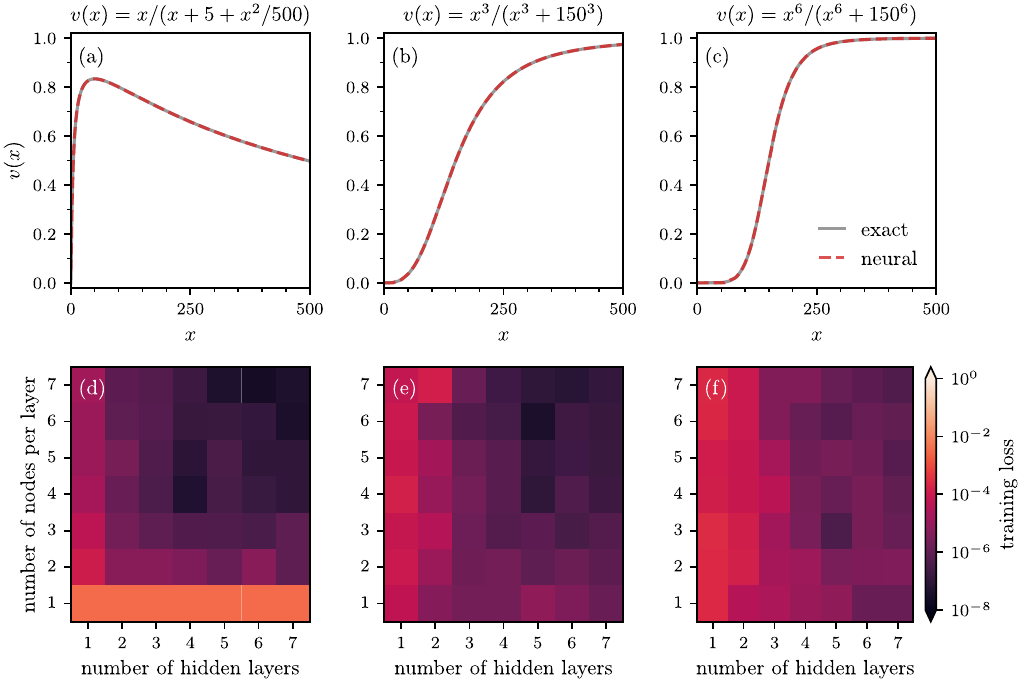}
    \caption{Approximation of three 1D target processes by $\mathrm{NNP}$s and the corresponding training losses across network architectures. (a--c) The solid grey lines denote the target processes $v(x)$ as functions of the state variable $x$, and the dashed red lines show the best-performing neural-network approximations over 100 random initializations of the network parameters. The parameters are initialized from a uniform distribution $\mathcal{U}(-d_{\rm in}^{-1/2},d_{\rm in}^{-1/2})$, where $d_{\rm in}$ is the number of input features. The target functions are displayed above the corresponding panels. (d--f) Heatmaps of the best training loss for each architecture as a function of the number of hidden layers and the number of nodes per layer. The color scale is logarithmic, with darker regions indicating lower loss.}
    \label{fig:nnp_approximation_1d}
\end{figure}

\begin{figure}
    \centering
    \includegraphics{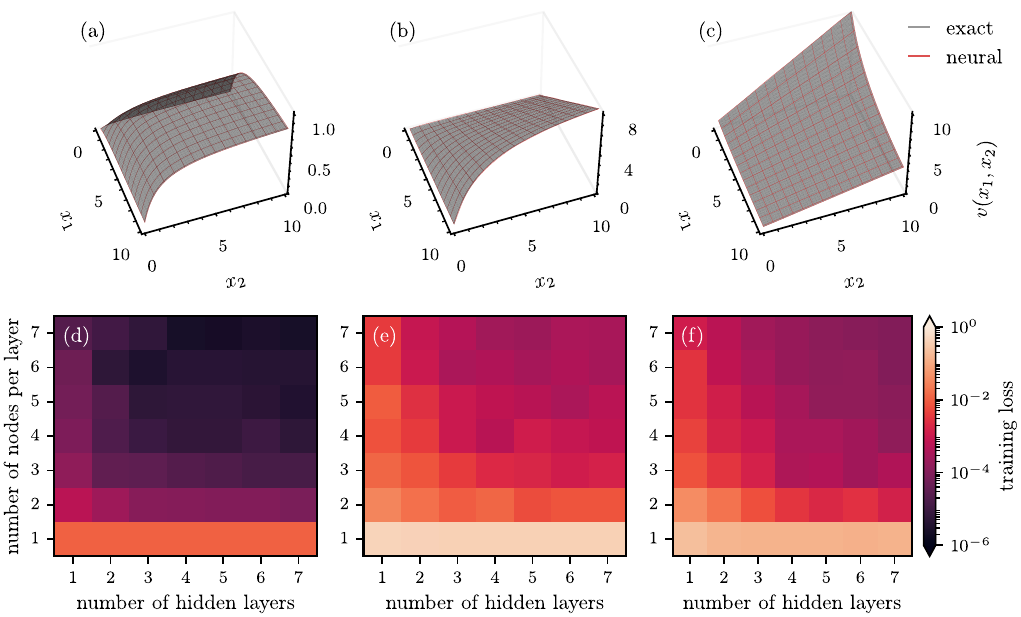}
    \caption{Approximation of three 2D target processes by $\mathrm{NNP}$s and the corresponding training losses across network architectures. (a--c) The grey surfaces denote the target processes: (a) $v(x_1,x_2)= x_1 x_2 / (1.35 + 0.9 x_1 + 1.5 x_2 + x_1 x_2)$, (b) $v(x_1,x_2)=x_1 x_2 / (2 + x_2)$, and (c) $v(x_1,x_2) = (1 + x_2)/(1 + 0.25 x_1)$, each as a function of the variables $x_1$ and $x_2$. The red wireframes show the best-performing neural-network approximations over 100 random initializations of the network parameters. The parameters are initialized from a uniform distribution $\mathcal{U}(-d_{\rm in}^{-1/2},d_{\rm in}^{-1/2})$, where $d_{\rm in}$ is the number of input features. (d--f) Heatmaps of the best training loss for each architecture as a function of the number of hidden layers and the number of nodes per layer. The color scale is logarithmic, with darker regions indicating lower loss.}
    \label{fig:nnp_approximation_2d}
\end{figure}

In this section, we examine how network depth and width affect the ability of an $\mathrm{NNP}$ to approximate 1D and 2D biological processes. To this end, we generated datasets with 1,000 samples for the 1D processes and 5,000 samples for the 2D processes. For each target process, we trained fully connected $\mathrm{NNP}$s with up to seven hidden layers and up to seven nodes per hidden layer. Each architecture was trained 100 times from different random initializations. The best-performing fits are shown in Figs.~\ref{fig:nnp_approximation_1d}(a--c) and~\ref{fig:nnp_approximation_2d}(a--c), while panels~(d)--(f) summarize the corresponding training losses across architectures.

In Fig.~\ref{fig:nnp_approximation_1d}, we consider three representative 1D processes: (a) a Haldane rate law and (b,c) Hill functions with Hill coefficients $3$ and $6$, respectively. For the Haldane process, we obtained low losses already for networks with three to four hidden layers and three to four nodes per layer. The Hill functions required larger architectures, and the required network size increased with steepness. The Hill function with coefficient $3$ was well approximated by networks with five hidden layers and four to six nodes per layer, whereas the steeper case with coefficient $6$ was not captured equally well within the $7\times7$ architecture grid.

Similarly, $\mathrm{NNP}$s also approximated well the three representative 2D processes shown in Fig.~\ref{fig:nnp_approximation_2d}: (a) a bisubstrate Michaelis--Menten rate law, (b) a Monod-type process, and (c) a gene-transcription process with one activator and one inhibitor based on the formulation proposed by Ac\'on et al.~\cite{acon2021myc}. Of these three examples, the bisubstrate Michaelis--Menten process was already well-approximated by architectures with at least three hidden layers and four nodes per layer. The Monod-type and gene-transcription examples required somewhat larger architectures to achieve a similar approximation error. For architectures with comparable approximation error, the computational cost increased more strongly with network depth than with width, suggesting that shallower but wider architectures may be more computationally efficient in this setting (see Appendix~\ref{app:runtime}). Overall, these numerical experiments indicate that relatively small-scale $\mathrm{NNP}$s are sufficient for several standard biological rate laws, whereas sharper nonlinearities may require larger networks.

\section{Biochemically informed neural ordinary differential equations}
With the $\mathrm{NNP}$ representation introduced in the previous section, we can now construct neural ODE models that preserve the process-based structure of stoichiometric systems. A biochemically informed neural ordinary differential equation (BINODE) represents the dynamics of a given system as a collection of nonnegative neural network processes whose outputs are mapped linearly to the state derivatives. In this way, individual biochemical processes remain explicit components of the model while being represented flexibly by neural networks.

We write a BINODE with $n$ state variables and $k$ processes as
\begin{equation}
    \frac{\mathrm{d}\hat{X}}{\mathrm{d}t}=W\,\mathcal{V}(\hat{X};\Theta)\,,
    \label{eq:BINODE}
\end{equation}
where $W\in\mathbb{R}^{n\times k}$ is a linear output layer analogous to a stoichiometric matrix, and $\Theta = (\theta_1,\dots,\theta_k)$ denotes the $\mathrm{NNP}$ parameters. The vector of process rates is
\begin{equation}
    \mathcal{V}(\hat{X};\Theta)=\big(\mathrm{NNP}_1(\hat{X};\theta_1),\dots,\mathrm{NNP}_k(\hat{X};\theta_k)\big)^\top \in \mathbb{R}^k\,.
\end{equation}
For notational simplicity, each process is written as a function of the full approximated state vector $\hat{X}\equiv \hat{X}(t) = (\hat{x}_1(t), \dots, \hat{x}_n(t))^\top \in \mathbb{R}^n$, although in practice an $\mathrm{NNP}$ may depend on only a subset of the state variables.

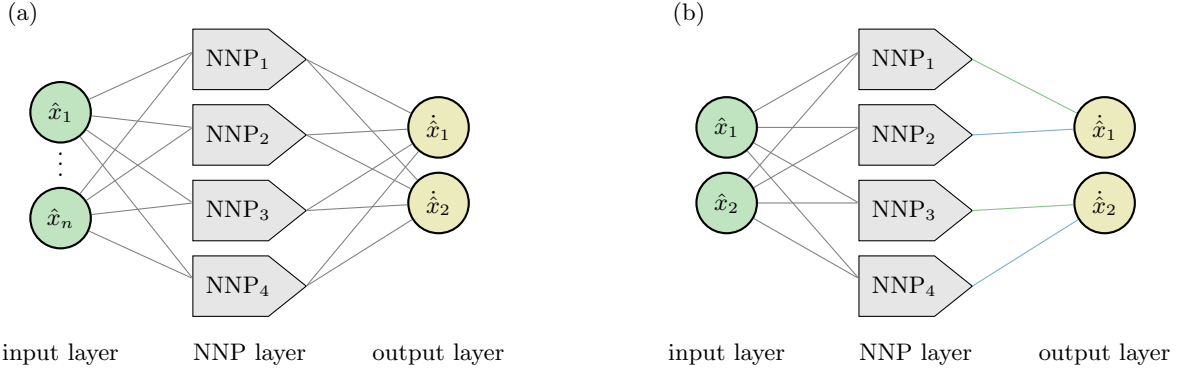
\begin{figure}
    \centering
    \input{BINODE.tikz}
    \hspace{5em}
    \input{BINODE_Ssystem.tikz}
    \caption{Schematics of BINODEs. (a) General BINODE architecture with state variables $x_1,\dots,x_n$ and a process layer composed of multiple $\mathrm{NNP}$s. The output of the network is the vector of state derivatives $\dot{x}_1,\dots,\dot{x}_n$. (b) BINODE with an S-system structure for two state variables. Each derivative is represented by the difference between two processes, one aggregating positive contributions and one aggregating negative contributions [see Eq.~\eqref{eq:Ssystem}]. Green and blue edges indicate positive and negative output-layer weights, respectively. Additional structural restrictions, such as monotonicity of the constituent processes, can be imposed to mimic power-law behavior. In both panels, the process layer is connected to the outputs through a linear stoichiometric layer with real-valued weights and no biases.}
    \label{fig:BINODE}
\end{figure}

At the system level, the BINODE can be viewed as a three-stage architecture [see Fig.~\ref{fig:BINODE}(a)]. The first stage contains the state variables $x_1,\dots,x_n$. The second stage consists of the constituent processes, each represented by an $\mathrm{NNP}$. The third stage is a linear output layer, without bias terms, that combines the process rates to produce $\dot{x}_1,\dots,\dot{x}_n$. 

Depending on the available prior knowledge, biochemical information can be incorporated by restricting which variables are supplied to each $\mathrm{NNP}$, fixing or constraining selected entries of $W$, or imposing structural constraints such as monotonicity on individual processes. For example, monotonicity in selected components can be supported by constraining neural network parameters to be nonnegative (\eg, via squaring the parameters). Additionally, input convex neural networks~\cite{DBLP:conf/icml/AmosXK17} provide a principled way to enforce convex dependence of the $\mathrm{NNP}$ output on a subset of inputs, and may be useful in applications where such structure is desired.

BINODEs retain the process-based organization of mechanistic ODE models while offering greater flexibility than closed-form rate laws. BINODEs are also more transparent than fully black-box neural ODEs for two reasons. First, the linear output layer plays the role of an effective stoichiometric matrix, while the $\mathrm{NNP}$s represent individual process rates. As a result, each learned process can be visualized and analyzed as a function $v(X)$, analogously to a mass-action term (\eg, $v(x_1,x_2)=x_1x_2$), a Michaelis--Menten rate law (\eg, $v(x_1)=V_{\rm max}x_1/(x_1+K_{\rm m})$), or any of the other process representations listed in Tables~\ref{tab:processes} and \ref{tab:processes_appendix}. Second, the dependence structure of each process can be examined directly. If a variable is omitted from the input set of a given $\mathrm{NNP}$, or if all of its input connections to that process are constrained to zero, then the corresponding process is independent of that variable. This is analogous to a zero kinetic order in the power-law model~\eqref{eq:power_law_process}.

As an illustrative special case, the S-system formalism from biochemical systems theory~\cite{savageau1969iochemical,savageau1969biochemical,savageau1970biochemical},
\begin{equation}
    \frac{\mathrm{d}x_i}{\mathrm{d}t}=\alpha_i \prod_j x_j^{g_{ij}}-\beta_i \prod_j x_j^{h_{ij}}\quad \mathrm{for}~i\in\{1,2,\dots,n\}\,,
    \label{eq:Ssystem}
\end{equation}
can be represented as a BINODE in which only a subset of the weights in the linear output layer are nonzero [see Fig.~\ref{fig:BINODE}(b)]. In an S-system, each state variable is modeled by two aggregated processes: a positive term collecting all incoming contributions and a negative term collecting all outgoing contributions. Both are represented by products of power-law functions. The corresponding neural S-system therefore uses two $\mathrm{NNP}$s per state variable, and each output $\dot{x}_i$ is given by the weighted difference of its associated positive and negative processes. We show a two-variable example in Fig.~\ref{fig:BINODE}(b), but the construction extends directly to higher-dimensional systems. Unlike the classical S-system, however, the neural version does not in general retain the algebraic steady-state properties that make the original formalism analytically convenient.

To train BINODEs, we form mini-batches by sampling multiple starting points from observed trajectories and integrating the model forward over a short time horizon of $H$ time steps. The training loss is computed as the mean squared error over these short rollout segments,
\begin{equation}
\mathcal{L}(\Theta, W)=
\frac{1}{B H}
\sum_{i=1}^B
\sum_{h=1}^H
\left\|
\hat{X}^{(i)}(t_i + h\Delta t)
-X^{(i)}(t_i + h\Delta t)
\right\|_2^2 \,,
\end{equation}
where $B$ is the batch size, $t_i$ are sampled starting times, and $\hat{X}^{(i)}$ denotes the model prediction obtained by integrating Eq.~\eqref{eq:BINODE} forward from the initial condition $X^{(i)}(t_i)$ over $H$ steps. We minimize $\mathcal{L}(\Theta, W)$ using the Adam optimizer with learning rates in the range $10^{-4}$ to $10^{-1}$.
\section{Applications in biological modeling}
In this section, we illustrate the BINODE framework on four representative systems arising in biological modeling. In the first two examples, we examine the ability of BINODEs to approximate Monod and Lotka--Volterra systems as paradigmatic models of biological dynamics. We then consider pharmacokinetics and ultradian endocrine dynamics as two complementary examples to connect our approach to the benchmarks considered in~\cite{AI-Aristotle} and to illustrate how BINODEs can be used in partially mechanistic settings, in which only selected process terms are replaced by neural representations.
\subsection{Bioreactor model}
\label{sec:monod}
As a first test case, we consider the bioreactor model
\begin{align}
\begin{split}
        \frac{\mathrm{d}x_1}{\mathrm{d}t} &= \mu(x_2) x_1 - k_{\rm d} x_1 \\
    \frac{\mathrm{d}x_2}{\mathrm{d}t} &= -\frac{1}{Y_{x_1 x_2}} \mu(x_2) x_1\,,
\end{split}
    \label{eq:bioreactor}
\end{align}
where $x_1 \equiv x_1(t)$ and $x_2 \equiv x_2(t)$ denote the biomass concentration and substrate concentration at time $t$, respectively. The parameter $k_{\rm d}$ is the biomass decay rate, $\mu(x_2)$ is the specific biomass growth rate, and $Y_{x_1x_2}$ is the biomass growth yield.

We model the specific growth rate by the Monod equation
\begin{equation}
    \mu(x_2)=\mu_{\rm max}\frac{x_2}{x_2+K_{x_2}}\,,
    \label{eq:monod}
\end{equation}
where $\mu_{\rm max}$ is the maximum specific growth rate and $K_{x_2}$ is the half-saturation constant. We list the parameter values used in our simulations in Table~\ref{tab:monod}.
\begin{table}
\footnotesize
\centering
\renewcommand*{\arraystretch}{1.6}
\begin{tabular}{>{\centering\arraybackslash} m{10em} >{\raggedright\arraybackslash} m{15em} >{\raggedright\arraybackslash} m{8em}}\toprule
\textbf{Symbol} & \textbf{Definition} & \textbf{Value}
\\[1pt] \midrule
\,\,\, $\mu_{\rm max}$\,\, & maximum specific growth rate & $0.86~\mathrm{h}^{-1}$ \\[1pt]
\,\,\, $K_{x_2}$\,\, & half-saturation constant & $0.0138~\mathrm{kg}/\mathrm{m}^3$ \\[1pt]
\,\,\, $Y_{x_1 x_2}$\,\, & biomass growth yield & $1.28$ \\[1pt]
\,\,\, $k_{\rm d}$\,\, & biomass decay rate & $3\times 10^{-2}~\mathrm{h}^{-1}$ \\[1pt] \bottomrule
\end{tabular}
\vspace{1mm}
\caption{Parameters used in the Monod bioreactor model [see Eqs.~\eqref{eq:bioreactor} and~\eqref{eq:monod}]. }
\label{tab:monod}
\end{table}

Introducing the nonnegative process rates
\begin{equation}
    v_1(x_1,x_2)=\mu(x_2)x_1
    \quad \text{and} \quad
    v_2(x_1)=k_{\rm d}x_1\,,
    \label{eq:monod_processes}
\end{equation}
we can write the system \eqref{eq:bioreactor} in stoichiometric form as
\begin{equation}
    \frac{\mathrm{d}}{\mathrm{d}t}
    \begin{pmatrix}
        x_1 \\
        x_2
    \end{pmatrix}
    =
    \begin{pmatrix}
        1 & -1 \\
        -Y_{x_1 x_2}^{-1} & 0
    \end{pmatrix}
    \begin{pmatrix}
        v_1(x_1,x_2) \\
        v_2(x_1)
    \end{pmatrix}\,.
    \label{eq:monod_matrix}
\end{equation}

\begin{figure}
    \centering
    \input{Monod_BINODE.tikz}
    \caption{Schematic of the BINODE used to learn the dynamics of the Monod model.}
    \label{fig:MonodBINODE}
\end{figure}
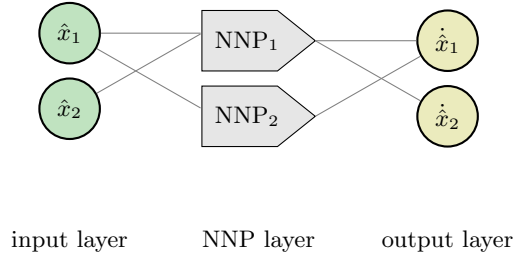

To learn this system with a BINODE, we employ two $\mathrm{NNP}$s and write
\begin{align}
 \frac{\mathrm{d}}{\mathrm{d}t}
\begin{pmatrix}
\hat{x}_1 \\
\hat{x}_2
\end{pmatrix}
&= \begin{pmatrix}
w_{11} & w_{12}  \\
w_{21} & 0
\end{pmatrix}
\begin{pmatrix}
\mathrm{NNP}_1\big(\hat{x}_1,\hat{x}_2;\theta_1\big) \\
\mathrm{NNP}_2\big(\hat{x}_1;\theta_2\big)
\end{pmatrix}\,,
    \label{eq:MonodBINODEeq}
\end{align}
where $\mathrm{NNP}_1$ captures biomass growth and substrate depletion, while $\mathrm{NNP}_2$ captures biomass decay. We implement each $\mathrm{NNP}$ as a fully connected network with 5 hidden layers of width 5, ELU activation functions, and a softplus output to ensure positive process rates. To encode the expected process structure, we use input masks to restrict the dependence of the processes on the state variables, so that the first process depends on $(\hat{x}_1,\hat{x}_2)$ and the second only on $\hat{x}_1$. A masked linear output layer maps the two processes to the two state derivatives. We show a schematic of this BINODE in Fig.~\ref{fig:MonodBINODE}.

\begin{figure}
    \centering
    \includegraphics{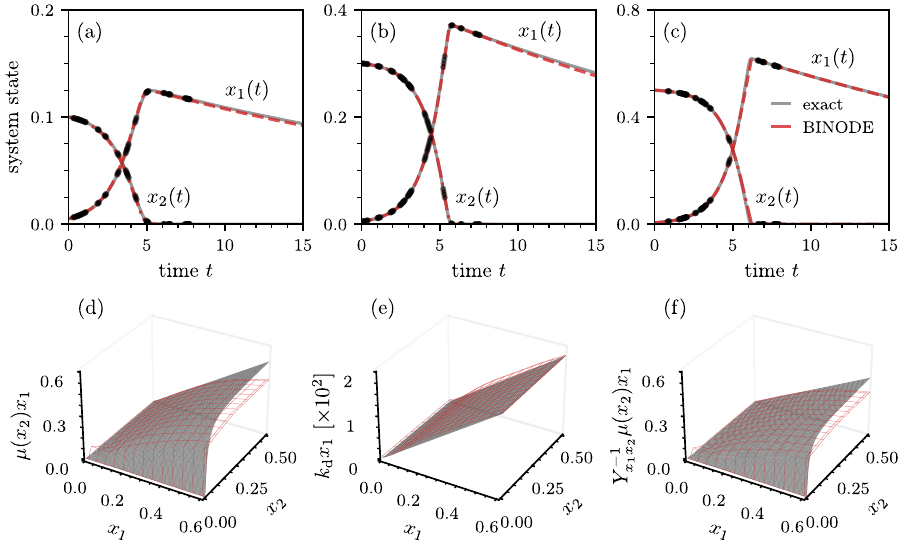}
    \caption{Comparison of BINODE predictions with the reference Monod model. (a--c) Time evolution of the state variables $x_1(t)$ and $x_2(t)$ for three different initial conditions, with units of $\mathrm{kg}/\mathrm{m}^3$. Black dots denote training observations, grey lines denote the reference trajectories, and red lines denote the BINODE trajectories. (d--f) Grey surfaces denote the target process contributions $v_1(x_1,x_2)=\mu(x_2)x_1$, $v_2(x_1)=k_{\rm d}x_1$, and $Y_{x_1 x_2}^{-1}v_1(x_1,x_2)=Y_{x_1 x_2}^{-1}\mu(x_2)x_1$, while the red wireframes show the corresponding learned contributions $w_{11}\,\mathrm{NNP}_1\big(\hat{x}_1,\hat{x}_2;\theta_1\big)$, $-w_{12}\,\mathrm{NNP}_2\big(\hat{x}_1;\theta_2\big)$, and $-w_{21}\,\mathrm{NNP}_1\big(\hat{x}_1,\hat{x}_2;\theta_1\big)$.}
    \label{fig:monod}
\end{figure}

The BINODE is trained on samples from trajectories with initial conditions $(0.005,0.1)$, $(0.005,0.3)$, and $(0.005,0.5)$, where both state variables are measured in $\mathrm{kg}/\mathrm{m}^3$ [see Fig.~\ref{fig:monod}(a--c)]. We use a training batch size of $B=20$ and a rollout horizon of $H=4$ time steps. As shown in Fig.~\ref{fig:monod}(a--c), the learned model reproduces the reference trajectories and approximates the corresponding process contributions well. Because an $\mathrm{NNP}$ output and its associated output-layer weight can be rescaled jointly without changing the resulting vector field, the comparison in Fig.~\ref{fig:monod}(d--f) is made at the level of the combined process contributions rather than the raw $\mathrm{NNP}$ outputs. In Appendix~\ref{app:monod_empirical}, we apply this BINODE setup to empirical biodegradation data to study process reconstruction in sparse and noisy settings.

To place these results in context, we compare the BINODE approach with SINDy~\cite{desilva2020,Kaptanoglu2022,SINDy2016,Kaptanoglu_pysindy}, a widely used method for data-driven model discovery. Specifically, we examine its ability to generate polynomial approximations for the bioreactor example. Using a single training trajectory generated from Eq.~\eqref{eq:bioreactor}, with the parameter values given in Table~\ref{tab:monod} and initial condition $(0.005, 0.7)$, the results are mixed.

On the one hand, SINDy is highly efficient at generating candidate models. On the other hand, not all identified models yield ODE systems that can be successfully simulated. In this example, the second-, fourth-, and sixth-order polynomial approximations fail to produce integrable models. The third- and fifth-order approximations can be simulated, with the fifth-order model providing a better fit. The seventh-order approximation yields negative values, while the eighth-order model is not integrable. The ninth- and tenth-order approximations both integrate successfully and provide good fits.

These results highlight the challenges of approximating rational functions with polynomials, which form the standard library in SINDy. While SINDy includes extensions for identifying implicit ODE models that can represent rational functions, these approaches also yield mixed results.
\subsection{Lotka--Volterra model}
We next consider the classical Lotka--Volterra predator--prey system,
\begin{align}
\begin{split}
    \frac{\mathrm{d}x_1}{\mathrm{d}t} &= \alpha x_1 - \beta x_1 x_2 \\
    \frac{\mathrm{d}x_2}{\mathrm{d}t} &= \gamma x_1 x_2 - \delta x_2\,,
\end{split}
    \label{eq:predator-prey}
\end{align}
where $x_1$ and $x_2$ denote the prey and predator populations, respectively, while $\alpha$, $\beta$, $\gamma$, and $\delta$ are the prey growth, predation, predator growth, and predator death rates.

Writing the processes as
\begin{equation}
    v_1(x_1)=x_1\,,
    \qquad
    v_2(x_2)=x_2\,,
    \qquad
    v_3(x_1,x_2)=x_1x_2\,,
    \label{eq:lv_processes}
\end{equation}
the Lotka--Volterra system in stoichiometric form is
\begin{equation}
    \frac{\mathrm{d}}{\mathrm{d}t}
    \begin{pmatrix}
        x_1 \\
        x_2
    \end{pmatrix}
    =
    \begin{pmatrix}
        \alpha & 0 & -\beta \\
        0 & -\delta & \gamma
    \end{pmatrix}
    \begin{pmatrix}
        v_1(x_1) \\
        v_2(x_2) \\
        v_3(x_1,x_2)
    \end{pmatrix}\,,
    \label{eq:lv_matrix}
\end{equation}
which separates prey growth, predator decay, and predation into three distinct processes.

To learn this system, we use a BINODE of the form
\begin{align}
\frac{\mathrm{d}}{\mathrm{d}t}
\begin{pmatrix}
    \hat{x}_1 \\
    \hat{x}_2
\end{pmatrix}
&=
\begin{pmatrix}
    w_{11} & 0 & w_{13} \\
    0 & w_{22} & w_{23}
\end{pmatrix}
\begin{pmatrix}
    \mathrm{NNP}_1\big(\hat{x}_1;\theta_1\big) \\
    \mathrm{NNP}_2\big(\hat{x}_2;\theta_2\big) \\
    \mathrm{NNP}_3\big(\hat{x}_1,\hat{x}_2;\theta_3\big)
\end{pmatrix}\,.
\label{eq:LVBINODEeq}
\end{align}
Here, $\mathrm{NNP}_1$, $\mathrm{NNP}_2$, and $\mathrm{NNP}_3$ are intended to capture prey growth, predator decay, and predation, respectively. We implement each $\mathrm{NNP}$ as a fully connected network with 5 hidden layers of width 5, ELU activations, and a softplus output to ensure positive outputs. Input masks restrict the modules to depend on $\hat{x}_1$, $\hat{x}_2$, or both, and a masked linear output layer maps the three processes to the two state derivatives. We show a schematic of this BINODE in Fig.~\ref{fig:LVBINODE}.

For this example, we generated three training trajectories from Eq.~\eqref{eq:predator-prey} with parameters $\alpha=\beta=\gamma=\delta=1$ and initial conditions $(1.6, 0.4)$, $(0.5, 1.5)$, and $(1.7, 1.7)$ [see Fig.~\ref{fig:lv_example}(a--c)]. These trajectories jointly sample a broad region of the phase space. We use a batch size of $B=40$ and a rollout horizon of $H=10$ time steps. As shown in Fig.~\ref{fig:lv_example}, the learned BINODE reproduces the oscillatory dynamics and recovers a process decomposition with one approximately linear function of $x_1$, one approximately linear function of $x_2$, and one bilinear interaction surface in $(x_1,x_2)$, consistent with the structure of the reference model.
\begin{figure}
    \centering
    \input{LV_BINODE.tikz}
    \caption{Schematic of the BINODE used to learn the dynamics of the Lotka--Volterra system.}
    \label{fig:LVBINODE}
\end{figure}
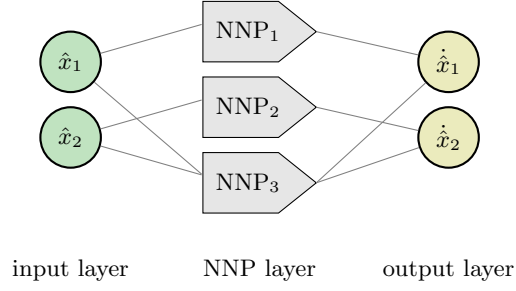
\begin{figure}
    \centering
    \includegraphics{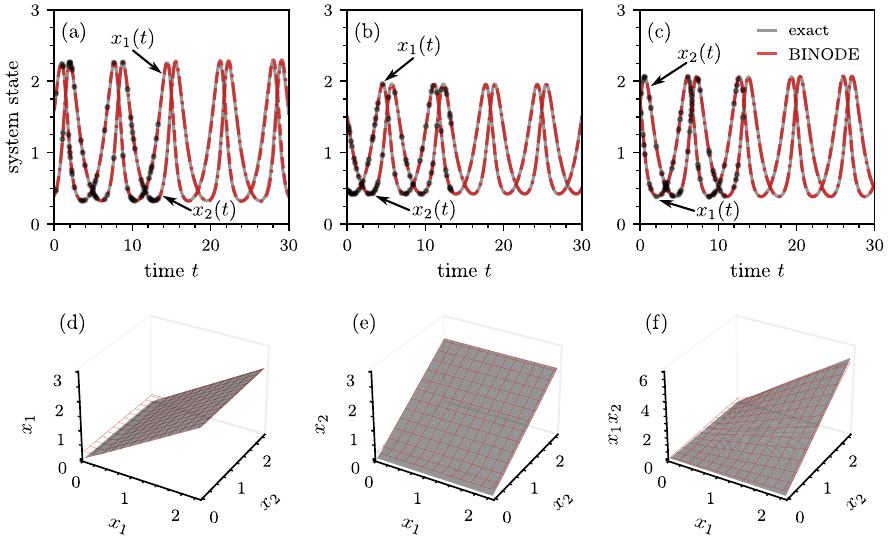}
    \caption{Comparison of BINODE predictions with the reference Lotka--Volterra model. (a--c) Time evolution of the state variables $x_1(t)$ and $x_2(t)$ for three different initial conditions. Black dots denote training observations; for visual clarity, we show only the first time point of each rollout segment. Grey lines denote the reference trajectories, and red lines denote the BINODE trajectories. (d--f) Grey surfaces denote the target process contributions $v_1(x_1)=x_1$, $v_2(x_2)=x_2$, and $v_3(x_1,x_2)=x_1x_2$, while the red wireframes show the corresponding learned contributions $w_{11}\,\mathrm{NNP}_1\big(\hat{x}_1;\theta_1\big)$, $-w_{22}\,\mathrm{NNP}_2\big(\hat{x}_2;\theta_2\big)$, and $w_{23}\,\mathrm{NNP}_3\big(\hat{x}_1,\hat{x}_2;\theta_3\big)$.}
    \label{fig:lv_example}
\end{figure}
\subsection{Pharmacokinetics model}
\begin{figure}
    \centering
    \includegraphics{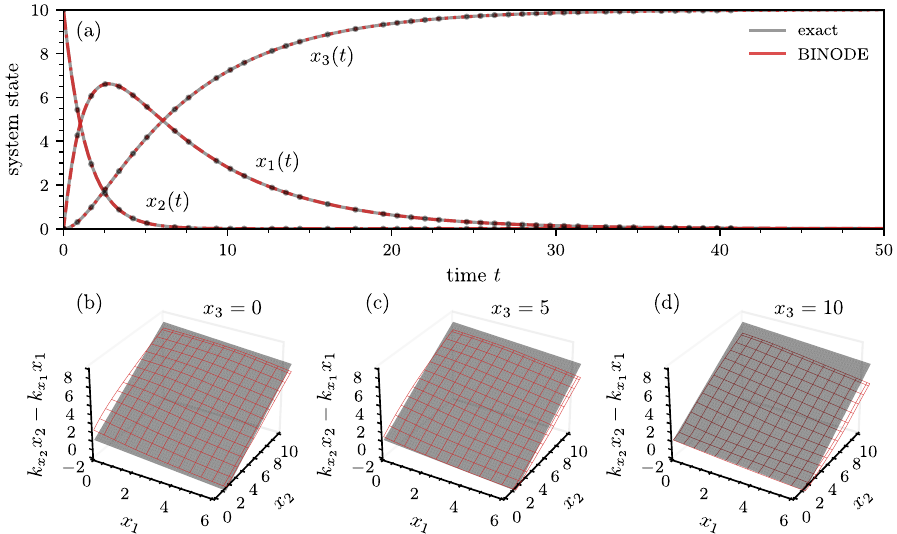}
    \caption{Comparison of BINODE predictions with the reference pharmacokinetics model. (a) Time evolution of the state variables $x_1(t)$, $x_2(t)$, and $x_3(t)$. Black dots denote training observations; for visual clarity, we show only the first time point of each rollout segment. Grey lines denote the reference trajectories, and red lines denote the BINODE trajectories. (b--d) Grey surfaces denote the target process term $k_{x_2}x_2-k_{x_1}x_1$, while the red wireframes show the learned rate of change $w\,\mathrm{NNP}(\hat{x}_1,\hat{x}_2,\hat{x}_3;\theta)$ evaluated at fixed values of $x_3=0$, $x_3=5$, and $x_3=10$, respectively.}
    \label{fig:pharmacokinetics}
\end{figure}
For comparison with AI-Aristotle~\cite{AI-Aristotle}, we consider the single-dose compartmental pharmacokinetics benchmark used in that work. AI-Aristotle combines the eXtreme Theory of Functional Connections, physics-informed neural networks, and symbolic regression to identify governing equations from data~\cite{AI-Aristotle,schiassi2021X-TFC,de2022PINNs,virgolin2022symbolic,fronk2023interpretable,gplearn}. The reference pharmacokinetics model~\cite{AI-Aristotle,barnes2011mathematical} is
\begin{align}
    \frac{\mathrm{d}x_1}{\mathrm{d}t} &= k_{x_2} x_2 - k_{x_1} x_1 \\
    \frac{\mathrm{d}x_2}{\mathrm{d}t} &= -k_{x_2} x_2 \\
    \frac{\mathrm{d}x_3}{\mathrm{d}t} &= k_{x_1} x_1\,,
    \label{eq:PKmodel}
\end{align}
with parameters $k_{x_1}=0.15~\mathrm{h}^{-1}$ and $k_{x_2}=0.72~\mathrm{h}^{-1}$. The states $x_1$, $x_2$, and $x_3$ represent the amounts of drug in the bloodstream, gastrointestinal tract, and urinary tract, respectively.

In the BINODE formulation, we replace the rate of change of $x_1$ with an $\mathrm{NNP}$ while keeping the remaining two equations fixed. That is,
\begin{align}
    \frac{\mathrm{d}\hat{x}_1}{\mathrm{d}t} &= w\,\mathrm{NNP}(\hat{x}_1,\hat{x}_2,\hat{x}_3;\theta) \\
    \frac{\mathrm{d}\hat{x}_2}{\mathrm{d}t} &= -k_{x_2} \hat{x}_2 \\
    \frac{\mathrm{d}\hat{x}_3}{\mathrm{d}t} &= k_{x_1} \hat{x}_1\,.
    \label{eq:learntPKmodel}
\end{align}
To match the representation used in~\cite{AI-Aristotle}, we model the rate of change of $x_1$ by a single $\mathrm{NNP}$, implemented as a fully connected network with 5 hidden layers of width 5, ELU activation functions, and a real-valued output. In contrast to~\cite{AI-Aristotle}, where the dynamics of $x_1$ are represented as an explicit function of time, we instead learn a state-dependent mapping. This yields an autonomous representation of the dynamics, in contrast to the explicitly time-dependent formulation in~\cite{AI-Aristotle}, and is therefore better suited to capturing the underlying mechanisms. Furthermore, to reflect incomplete prior knowledge of the governing equations, we allow dependence on all state variables $x_1$, $x_2$, and $x_3$. This tests whether the method identifies the correct dependence on $x_1$ and $x_2$, while suppressing spurious dependence on $x_3$.

In accordance with~\cite{AI-Aristotle}, we generate a single training trajectory on the time interval $[0,10]$ from the initial condition $(x_1(0),x_2(0),x_3(0))=(0,0.1~\mu\mathrm{g},0)$. We use a batch size of $B=40$ and a rollout horizon of $H=10$ time steps.

Figure~\ref{fig:pharmacokinetics}(a) shows good agreement between the reference model and the BINODE solution. Panels (b--d) indicate that the learned neural process recovers the linear dependence on $x_1$ and $x_2$, with very little dependence on $x_3$, consistent with the reference model~\eqref{eq:PKmodel}.
\subsection{Ultradian endocrine model}
\label{sec:ultradian}
\begin{figure}
    \centering
    \includegraphics{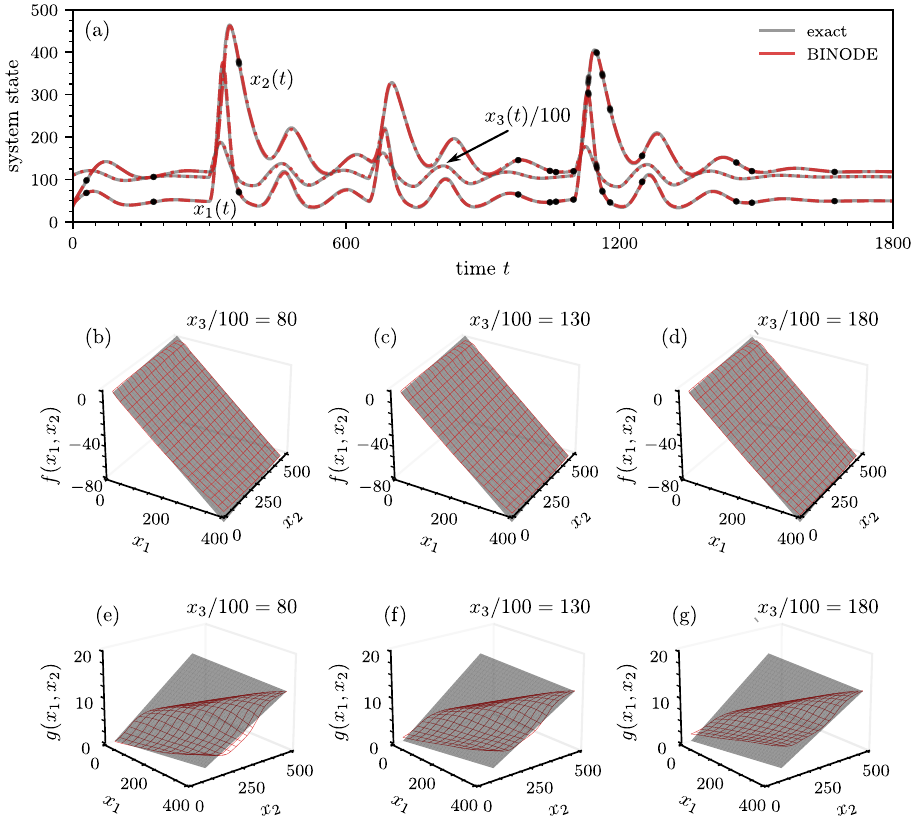}
    \caption{Comparison of BINODE predictions with the reference ultradian endocrine model. (a) Time evolution of the state variables $x_1(t)$, $x_2(t)$, and $x_3(t)/100$. The units are $\mu\mathrm{U}/\mathrm{mL}$, $\mu\mathrm{U}/\mathrm{mL}$, and $\mathrm{mg}$, respectively. The variable $x_3$ is rescaled by a factor of 100 for visualization. Grey lines denote the reference trajectories and red lines denote the BINODE trajectories. (b--d) Grey surfaces denote the first nonlinear target process term considered in the BINODE formulation, while the red wireframes show the corresponding learned neural process $\mathrm{NNP}_1(\hat{x}_1,\hat{x}_2,\hat{x}_3;\theta_1)$ at three fixed values of $x_3$. (e--g) Analogous comparison for the second nonlinear target process term and $\mathrm{NNP}_2(\hat{x}_1,\hat{x}_2,\hat{x}_3;\theta_2)$.}
    \label{fig:ultradian}
\end{figure}

Finally, we consider the six-dimensional ultradian endocrine model used in~\cite{AI-Aristotle}. It describes rhythmic oscillations in glucose and insulin levels on timescales typically shorter than 24 hours. This example combines transport terms, delayed feedback, and nonlinear physiological response functions. The dynamics are given by
\begin{align}
\begin{split}
\frac{\mathrm{d}x_1}{\mathrm{d}t} &= f_1(x_3) + f(x_1,x_2) \\
\frac{\mathrm{d}x_2}{\mathrm{d}t} &= g(x_1,x_2) \\
\frac{\mathrm{d}x_3}{\mathrm{d}t} &= f_4(x_6) + I_{\rm G}(t) - f_2(x_3) - f_3(x_2)\,x_3 \\
\frac{\mathrm{d}x_4}{\mathrm{d}t} &= \frac{x_1 - x_4}{t_{\rm d}} \\
\frac{\mathrm{d}x_5}{\mathrm{d}t} &= \frac{x_4 - x_5}{t_{\rm d}} \\
\frac{\mathrm{d}x_6}{\mathrm{d}t} &= \frac{x_5 - x_6}{t_{\rm d}}\,,
\end{split}
\end{align}
together with
\begin{align}
    \begin{split}
        f(x_1,x_2)&=- E\left(\frac{x_1}{V_1} - \frac{x_2}{V_2}\right) - \frac{x_1}{t_1}\\
        g(x_1,x_2)&=E\left(\frac{x_1}{V_1} - \frac{x_2}{V_2}\right) - \frac{x_2}{t_2}\,,
    \end{split}
\end{align}
\begin{align}
\begin{split}
        f_1(x_3) &= \frac{R_{\rm m}}{1 + \exp\!\left(-\frac{x_3}{V_3 C_1} + a_1\right)}\\
    f_2(x_3) &= U_{\rm b}\left(1 - \exp\!\left(-\frac{x_3}{C_2 V_3}\right)\right)  \\
    f_3(x_2) &= \frac{1}{C_3 V_3}
\left(
U_0 + \frac{U_m}{1 + (\kappa x_2)^{-\beta}}
\right) \\
   f_4(x_6) &= \frac{R_{\rm g}}{1 + \exp\!\left(\alpha\left(\frac{x_6}{C_5 V_1} - 1\right)\right)}\,,    
\end{split}
\end{align}
where $\kappa = \frac{1}{C_4}\left(\frac{1}{V_2} + \frac{1}{E t_2}\right)$. The three main state variables are plasma insulin $x_1(t)$, interstitial insulin $x_2(t)$, and glucose concentration $x_3(t)$. In addition, the state variables $x_4(t)$, $x_5(t)$, and $x_6(t)$ capture the delayed effect of insulin on glucose production~\cite{sturis1991computer}.

The functions $f_1(x_3)$, $f_2(x_3)$, $f_3(x_2)$, and $f_4(x_6)$ represent insulin secretion, insulin-independent glucose utilization, and two distinct components of insulin-dependent glucose utilization, respectively.

The nutritional input term is
\begin{equation}
    I_{\rm G}(t) = \sum_{t_i \le t} q_i\, k\, e^{-k (t - t_i)}\,,
\end{equation}
where $(t_i, q_i) = (300, 60)$, $(650, 40)$, $(1100, 50)$ $(\mathrm{min}, \mathrm{g})$. We list the remaining model parameters in Table~\ref{tab:ultradian_parameters} of Appendix~\ref{app:ultradian}. The initial conditions are $x_1(0)=36~\mu\mathrm{U}/\mathrm{mL}$, $x_2(0)=44~\mu\mathrm{U}/\mathrm{mL}$, $x_3(0)=11000~\mathrm{mg}$, $x_4(0)=0$, $x_5(0)=0$, and $x_6(0)=0$. The value $x_3(0)$ corresponds to an initial glucose concentration of $110~\mathrm{mg}/\mathrm{dL}$.

The learning task formulated in~\cite{AI-Aristotle} focuses on identifying time-dependent parameterizations of $f(\cdot)$ and $g(\cdot)$. Here, instead, we learn state-dependent mappings, yielding an autonomous representation of the dynamics. As in the prior section, we allow the $\mathrm{NNP}$s to depend on all state variables $x_1$, $x_2$, and $x_3$, enabling the model to identify the correct dependencies while suppressing spurious ones. We thus write
\begin{equation}
\begin{split}
    \frac{\mathrm{d}\hat{x}_1}{\mathrm{d}t} &= f_1(\hat{x}_3) + \mathrm{NNP}_1(\hat{x}_1,\hat{x}_2,\hat{x}_3;\theta_1) \\
    \frac{\mathrm{d}\hat{x}_2}{\mathrm{d}t} &= \mathrm{NNP}_2(\hat{x}_1,\hat{x}_2,\hat{x}_3;\theta_2)\,. \\
\end{split}
\end{equation}
The remaining mechanistic structure is kept fixed. To train the BINODE, we initially use a batch size of $B=5$ and a rollout horizon of $H=5$, adding further batches of the same size and horizon after 100 and approximately 1200 epochs. Figure~\ref{fig:ultradian}(a) shows that the trained BINODE reproduces the reference trajectories. Panels (b--d) indicate that $\mathrm{NNP}_1$ accurately approximates the surface $f(x_1,x_2)$ with negligible dependence on $x_3$, while panels (e--g) show that $\mathrm{NNP}_2$ approximates $g(x_1,x_2)$ well in regions where $x_1 \approx x_2$, but less accurately elsewhere, consistent with the similar behavior of $x_1$ and $x_2$ in the reference trajectories (\ie, the training data).
\section{Discussion and conclusions}
Dynamical systems in biochemistry and related fields are often written in stoichiometric form, involving a stoichiometric matrix and a corresponding process vector [see Eq.~\eqref{eq:NP}]. Depending on the application, a wide range of process representations has been developed, as summarized in Tables~\ref{tab:processes} and~\ref{tab:processes_appendix}, spanning more than a century since the foundational work of Waage and Guldberg~\cite{waage1864studier} on mass-action kinetics.

This variety of rate laws, together with the well-established approximation capabilities of artificial neural networks, motivates the use of neural network processes ($\mathrm{NNP}$s) as flexible building blocks for stoichiometric models. These $\mathrm{NNP}$s can be embedded in stoichiometric representations and combined to form biochemically informed neural ordinary differential equations (BINODEs). Using modern machine learning frameworks such as \texttt{PyTorch}, such models can be implemented efficiently via automatic differentiation.

After examining the approximation properties of $\mathrm{NNP}$s across different architectures, we applied BINODEs to bioreactor, Lotka--Volterra, pharmacokinetics, and ultradian endocrine models. In all cases, $\mathrm{NNP}$s with relatively small architectures, consisting of only a few layers and neurons per layer, were sufficient to approximate both the system dynamics and the underlying process surfaces. For the bioreactor example, we further demonstrated that BINODEs can be applied to empirical data, where they not only reproduce the observed trajectories but also yield interpretable process representations.

Given their flexibility and interpretability, the proposed BINODEs may be useful as surrogate models for complex biomedical systems~\cite{fonseca2025optimal}. For practical applications, it would be valuable to develop a software framework comprising pretrained $\mathrm{NNP}$ modules that can be combined and adapted by users to approximate empirical process data in a modular and interpretable manner.

Another area for future exploration is the use of the NNP representation for modeling the effects of other variables, like temperature, ionic strength, or light, on biological processes. The effect of temperature on enzyme activity is usually modeled using the Arrhenius equation~\cite{fonseca2011complex}, but its effects on other biological processes may not be well characterized. An NNP can be trained with independent variables for which data has been gathered, and the NNP can then be used to visualize the learned dependency.
\subsection*{Data availability statement}
Our source codes are publicly available at \url{https://gitlab.com/ComputationalScience/binode}.
\clearpage

\acknowledgements{LB acknowledges support from the hessian.AI Service Center (which is funded by the Federal Ministry of Research, Technology and Space, BMFTR; grant number 16IS22091) and the hessian.AI Innovation Lab (which is funded by the Hessian Ministry for Digital Strategy and Innovation; grant number S-DIW04/0013/003). LLF and RL acknowledge financial support from the Defense Advanced Research Projects Agency (grant HR00112220038), and the National Institutes of Health (grants R01 GM127909 and R01 AI135128). RL also acknowledges financial support from the National Institutes of Health (grant R01 HL169974).}

\newpage
\appendix
\section{Common biological processes}
\label{app:common_processes}
\begin{sidewaystable*}[p]
\centering
\small
\renewcommand{\arraystretch}{1.25}
\setlength{\tabcolsep}{5pt}

\begin{threeparttable}
\caption{Commonly used process representations in biochemistry and biology (continued).}
\label{tab:processes_appendix}

\begin{tabularx}{\textheight}{@{}
    >{\RaggedRight\arraybackslash}p{4.1cm}
    >{\centering\arraybackslash}m{8.4cm}
    >{\centering\arraybackslash}p{1.8cm}
    >{\RaggedRight\arraybackslash}X
@{}}
\toprule
\textbf{Process type} & \textbf{Mathematical form} & \textbf{Symmetric} & \textbf{Additional comments} \\
\midrule

Mass-action
& $\displaystyle \alpha \prod_i x_i^{n_i}$
& yes
& $\alpha$ is the rate constant, and $n_i$ is the stoichiometric coefficient of species $i$~\cite{waage1864studier}. \\

Competitive product inhibition
& $\displaystyle \frac{V_{\max}S}{S+K_{\mathrm m}\left(1+\frac{P}{K_{\mathrm m}^P}\right)}$
& no
& $V_{\max}$ is the maximum reaction rate, $K_{\mathrm m}$ is the Michaelis--Menten constant for the substrate, and $K_{\mathrm m}^P$ the product concentration at which the forward reaction rate is half-maximal~\cite{henri1903lois,henri2006theorie,cornish2013origins}. \\

\makecell[l]{Substrate inhibition \\ (Haldane equation)}
& $\displaystyle \frac{V_{\max}S}{S\left(1+\frac{S}{K_{\rm i}}\right)+K_{\mathrm m}}$
& n/a
& $V_{\max}$ is the maximum reaction rate, $K_{\mathrm m}$ is the Michaelis--Menten constant for the substrate, and $K_{\rm i}$ is the substrate inhibition constant~\cite{haldane1930enzymes}. \\

Hill
& $\displaystyle \frac{V_{\max}S^h}{S^h+K_{\mathrm m}^h}$
& n/a
& $V_{\max}$ is the maximum reaction rate, $K_{\mathrm m}$ is the Michaelis--Menten constant for substrate $S$, and $h$ is the Hill coefficient representing the cooperativity of the binding sites. \\

Reversible Hill
& $\displaystyle \frac{V_{\max}\frac{S}{K_{\mathrm m}^S}\left(\frac{S}{K_{\mathrm m}^S}+\frac{P}{K_{\mathrm m}^P}\right)^{h-1}}{1+\left(\frac{S}{K_{\mathrm m}^S}+\frac{P}{K_{\mathrm m}^P}\right)^h}$
& no
& $V_{\max}$ is the maximum reaction rate, $K_{\mathrm m}^S$ is the Michaelis--Menten constant for substrate $S$, $K_{\mathrm m}^P$ is the Michaelis--Menten constant for product $P$, and $h$ is the Hill coefficient representing the cooperativity of the binding sites~\cite{hofmeyr1997reversible}. \\

Saturable and cooperative formalism
& $\displaystyle V\prod_i \frac{x_i^{n_i}}{K_i+x_i^{n_i}}$
& yes
& $V$ represents the saturated reaction rate for large values of all $x_i$, and $K_i$ and $n_i$ are parameters associated with variable $x_i$, encoding sensitivity to and saturation with respect to $x_i$~\cite{sorribas2007cooperativity}. \\

Monod
& $\displaystyle \mu_{\max}\frac{S}{S+K_{\mathrm m}}\,X$
& no
& $\mu_{\max}$ is the maximum specific growth rate under saturating substrate concentrations, $K_{\mathrm m}$ is the Michaelis--Menten constant for substrate $S$, and $X$ is the biomass concentration~\cite{monod1958recherches,monod1950technique}. \\

Holling type I
& $\displaystyle aT_{\rm s}XP$
& no
& $a$ is a proportionality constant, $T_{\rm s}$ is the prey-searching time, $X$ is the prey density, and $P$ is the predator density. \\

Holling type II
& $\displaystyle \frac{aT_{\rm t}X}{1+abX}\,P$
& no
& $a$ is a proportionality constant, $T_{\rm t}$ is the total time available for searching and handling prey, $b$ is the handling time per prey item, $X$ is the prey density, and $P$ is the predator density~\cite{cs1959components,holling1959some}. \\

\bottomrule
\end{tabularx}
\end{threeparttable}
\end{sidewaystable*}
Table~\ref{tab:processes_appendix} continues the summary of commonly used process representations in biochemistry and biology.
\section{Runtime for different neural network sizes}
\label{app:runtime}
\begin{figure}
    \centering
    \includegraphics{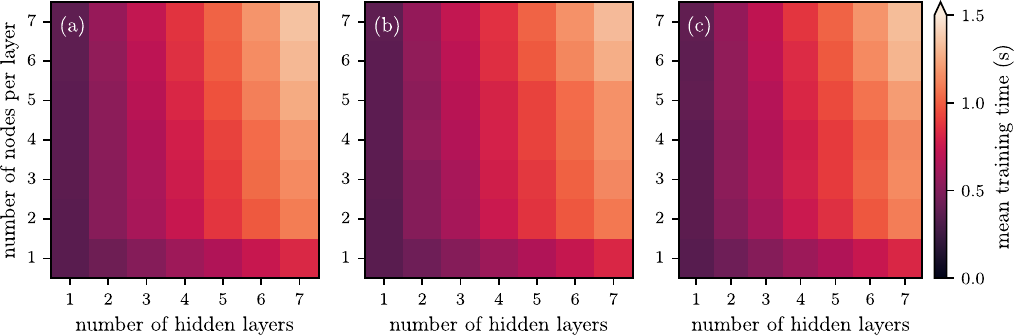}
    \caption{Mean training time for neural network models with varying architectures used to approximate three 1D target processes. (a–c) Heatmaps show the average training time (in seconds) across 100 random initializations for networks trained to approximate three 1D target processes (see Fig.~\ref{fig:nnp_approximation_1d}). The number of hidden layers and nodes per layer are varied along the vertical and horizontal axes, respectively.}
    \label{fig:MM_dynamic_1d_time}
\end{figure}
\begin{figure}
    \centering
    \includegraphics{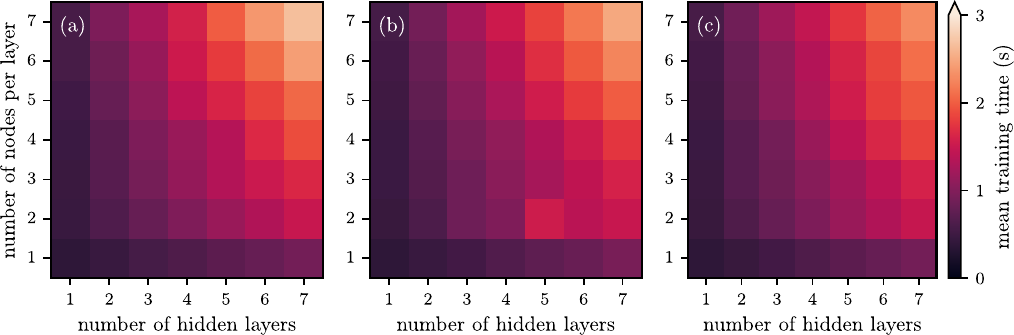}
    \caption{Mean training time for neural network models with varying architectures used to approximate three 2D target processes. (a–c) Heatmaps show the average training time (in seconds) across 100 random initializations for networks trained to approximate three 2D target processes (see Fig.~\ref{fig:nnp_approximation_2d}). The number of hidden layers and nodes per layer are varied along the vertical and horizontal axes, respectively.}
    \label{fig:MM_dynamic_2d_time}
\end{figure}
In Figs.~\ref{fig:MM_dynamic_1d_time} and \ref{fig:MM_dynamic_2d_time}, we show heatmaps of the mean runtime for the 1D and 2D target processes that we examined in Sec.~\ref{sec:characterizing_nnps}. In each heatmap, the number of hidden layers and the number of nodes per layer are varied along the horizontal and vertical axes, respectively.
\section{Application to empirical biodegradation data}
\label{app:monod_empirical}
\begin{figure}
    \centering
    \includegraphics{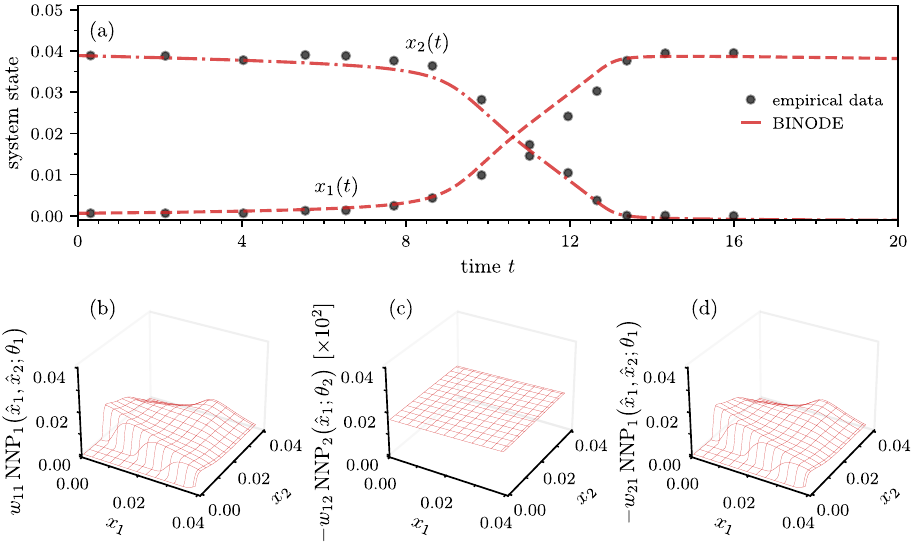}
    \caption{Application of the BINODE to empirical biodegradation data. (a) Time evolution of the state variables $x_1(t)$ (biomass) and $x_2(t)$ (toluene), with units of $\mathrm{kg}/\mathrm{m}^3$. Black dots indicate empirical observations, and red lines show the corresponding BINODE trajectories. (b--d) Red wireframes show the learned process contributions $w_{11}\,\mathrm{NNP}_1\big(\hat{x}_1,\hat{x}_2;\theta_1\big)$, $-w_{12}\,\mathrm{NNP}_2\big(\hat{x}_1;\theta_2\big)$, and $-w_{21}\,\mathrm{NNP}_1\big(\hat{x}_1,\hat{x}_2;\theta_1\big)$.}
    \label{fig:monod_empirical}
\end{figure}
To study process reconstruction in sparse and noisy settings, we apply the bioreactor BINODE from Section~\ref{sec:monod} to empirical biodegradation data from a toluene batch experiment with \textit{P.~putida} F1~\cite{reardon2000biodegradation}. We initialize the model using the bioreactor BINODE pre-trained on synthetic data and further train it on the empirical data, augmented with $B=15$ batches of rollout length $H=4$ generated from interpolated measurements.

As shown in Fig.~\ref{fig:monod_empirical}(a), the BINODE adapts to the observed dynamics, capturing the transient dynamics while indicating negligible decay. While the learned process surfaces in panels (b--d) do not resemble the canonical Monod form, they remain interpretable and capture key qualitative features of the underlying dynamics. They exhibit an initial increase followed by a decrease along the observed trajectory, consistent with the regime in which biomass ($x_1$) is initially low and substrate ($x_2$) is high. As biomass increases and substrate is depleted, the inferred growth rate increases and then declines. In addition, the learned decay contribution is close to zero, consistent with the observed data.

These results indicate that the BINODE is able to extract meaningful process structure even in the presence of sparse and noisy measurements.
\section{Ultradian endocrine model parameter values}
\label{app:ultradian}
\begin{table}[t]
\centering
\small
\renewcommand{\arraystretch}{1.2}
\setlength{\tabcolsep}{6pt}
\begin{tabular}{ccc}
\toprule
\textbf{Parameter} & \textbf{Nominal value} & \textbf{Unit} \\
\midrule
$V_1$   & 3     & $\mathrm{L}$ \\
$V_2$   & 11    & $\mathrm{L}$ \\
$V_3$   & 10    & $\mathrm{L}$ \\
$E$     & 0.2   & $\mathrm{L}\,\mathrm{min}^{-1}$ \\
$t_1$   & 6     & $\mathrm{min}$ \\
$t_2$   & 100   & $\mathrm{min}$ \\
$t_d$   & 12    & $\mathrm{min}$ \\
$k$     & 0.0083 & $\mathrm{min}^{-1}$ \\
$R_m$   & 209   & $\mathrm{mU}\,\mathrm{min}^{-1}$ \\
$a_1$   & 6.6   & -- \\
$C_1$   & 300   & $\mathrm{mg}\,\mathrm{L}^{-1}$ \\
$C_2$   & 144   & $\mathrm{mg}\,\mathrm{L}^{-1}$ \\
$C_3$   & 100   & $\mathrm{mg}\,\mathrm{L}^{-1}$ \\
$C_4$   & 80    & $\mathrm{mU}\,\mathrm{L}^{-1}$ \\
$C_5$   & 26    & $\mathrm{mU}\,\mathrm{L}^{-1}$ \\
$U_b$   & 72    & $\mathrm{mg}\,\mathrm{min}^{-1}$ \\
$U_0$   & 4     & $\mathrm{mg}\,\mathrm{min}^{-1}$ \\
$U_m$   & 90    & $\mathrm{mg}\,\mathrm{min}^{-1}$ \\
$R_g$   & 180   & $\mathrm{mg}\,\mathrm{min}^{-1}$ \\
$\alpha$ & 7.5  & -- \\
$\beta$  & 1.772 & -- \\
\bottomrule
\end{tabular}
\caption{Parameter values for the ultradian endocrine model.}
\label{tab:ultradian_parameters}
\end{table}
Table~\ref{tab:ultradian_parameters} lists the parameters of the ultradian endocrine model examined in Sec.~\ref{sec:ultradian}.
\clearpage
% Bibliography
\bibliographystyle{ieeetr}
\bibliography{refs.bib}
%%%%%%%%%%%%%%%%%%%%%%%%%%%%%%%%%%%%%%%%%%%%%%%%%%%%%%%%%%%%%%%%%%%%
%%%%%%%%%%%%%%%%%%%%%%%%%%%%%%%%%%%%%%%%%%%%%%%%%%%%%%%%%%%%%%%%%%%%
%%%%%%%%%%%%%%%%%%%%%%%%%%%%%%%%%%%%%%%%%%%%%%%%%%%%%%%%%%%%%%%%%%%%
\end{document}

%% file: nnp.tikz
\begin{tikzpicture}[
  neuron/.style={circle, draw=black, fill=#1, minimum size=8mm},
  input neuron/.style={neuron=tabgreen!30},
  hidden neuron/.style={neuron=tabblue!30},
  output neuron/.style={neuron=tabbrown!30},
  layer/.style={draw=none, rectangle}
]

\newcommand{\hidden}[3]{%
\node[inner sep=0pt, outer sep=0pt] (#3) at (#1,#2) {};
  \draw[fill=grey!20, draw=black]
    (#3)
    -- ++(0,1.5)
    -- ++(5,0)
    -- ++(2.5,-2)
    -- ++(-2.5,-2)
    -- ++(-5,0)
    -- cycle;
}
\hidden{1.35}{0}{H1}

% Input Layer
\node[input neuron] (I1) at (0,0.2) {$x_1$};
\node at (0,-0.38) {$\vdots$};
\node[input neuron] (I2) at (0,-1.2) {$x_n$};
\node[layer] at (-1.75,-0.5) (inputlabel) {input variables};

% Hidden Layers
\foreach \m [count=\y] in {1,2,3,4}
  \node[hidden neuron] (H1\m) at (2,2-\y) {$\sigma$};
\foreach \m [count=\y] in {1,2,3,4}
  \node[hidden neuron] (H2\m) at (4,2-\y) {$\sigma$};
\foreach \m [count=\y] in {1,2,3,4}
  \node[hidden neuron] (H3\m) at (6,2-\y) {$\sigma$};

% Output Layer
\node[output neuron] (O) at (8,-0.5) {$v$};
\node[layer, right=0.5cm of O] (outputlabel) {NNP output};

% Connections
\foreach \i in {1,2}
  \foreach \j in {1,2,3,4}
    \draw[gray] (I\i) -- (H1\j);

\foreach \i in {1,2,3,4}
  \foreach \j in {1,2,3,4}
    \draw[gray] (H1\i) -- (H2\j);

\foreach \i in {1,2,3,4}
  \foreach \j in {1,2,3,4}
    \draw[gray] (H2\i) -- (H3\j);

\foreach \i in {1,2,3,4}
  \draw[gray] (H3\i) -- (O);

% NNP label box
%\draw[rounded corners=10pt, thick, fill=black!20, opacity=0.3, draw=black] (1.0,1.75) rectangle (7.0,-2.75);

\node at (4.0,2.0) {neural network process (NNP)};

\end{tikzpicture}

%% file: BiNODE.tikz
\begin{tikzpicture}[node distance=1.5cm and 2.5cm, on grid, auto]

% Styles
\tikzstyle{input} = [circle, draw=black, fill=tabgreen!30, thick, minimum size=8mm]
\tikzstyle{output} = [circle, draw=black, fill=tabolive!30, thick, minimum size=8mm]
%\tikzstyle{hidden} = [rectangle, text=black, text opacity=1, draw=black, fill=black!40, fill opacity=0.3, thick, minimum size=8mm, rounded corners=3pt]
\newcommand{\hidden}[4]{% #1=x, #2=y, #3=name, #4=label
  \node[inner sep=0pt, outer sep=0pt] (#3) at (#1,#2) {};
  % Draw the polygon
  \draw[fill=grey!20, draw=black]
    (#3)
    -- ++(0,0.2*1.5)
    -- ++(0.2*5,0)
    -- ++(0.2*2.5,-0.2*2)
    -- ++(-0.2*2.5,-0.2*2)
    -- ++(-0.2*5,0)
    -- cycle;
  % Define right tip coordinate manually
  \coordinate (#3-right) at ({#1 + 1.5}, {#2 - 0.1});
  % Label
  \node at ({#1 + 0.6}, {#2 - 0.1}) {\textbf{#4}};
}

% Nodes
\node[input] (I1) at (0,0.2) {$\hat{x}_1$};
\node at (0,-0.38) {$\vdots$};
\node[input] (I2) at (0,-1.2) {$\hat{x}_n$};

%\node[hidden] (H1) at (2.5,1) {$\mathrm{NNP}_1$};
%\node[hidden] (H2) at (2.5,0) {$\mathrm{NNP}_2$};
%\node[hidden] (H3) at (2.5,-1) {$\mathrm{NNP}_3$};
%\node[hidden] (H4) at (2.5,-2) {$\mathrm{NNP}_4$};

\hidden{1.75}{1.0}{H1}{$\mathrm{NNP}_1$};
\hidden{1.75}{0.0}{H2}{$\mathrm{NNP}_2$};
\hidden{1.75}{-1.0}{H3}{$\mathrm{NNP}_3$};
\hidden{1.75}{-2.0}{H4}{$\mathrm{NNP}_4$};

\node[output] (O1) at (5,0) {$\dot{\hat{x}}_1$};
\node[output] (O2) at (5,-1) {$\dot{\hat{x}}_2$};

% Connections from Input to Hidden
\foreach \i in {I1,I2} {
  \foreach \h in {H1,H2,H3,H4} {
    \draw[gray] (\i) -- (\h);
  }
}

% Connections from Hidden to Output
\foreach \h in {H1,H2,H3,H4} {
  \foreach \o in {O1,O2} {
    \draw[gray] (\h-right) -- (\o);
  }
}

\node at (0,-3) (inputlabel) {input layer};
\node at (2.5,-3) (hiddenlabel) {NNP layer};
\node at (5,-3) (outputlabel) {output layer};

\node at (-0.5,1.5) {(a)};

\end{tikzpicture}

%% file: BiNODE_Ssystem.tikz
\begin{tikzpicture}[node distance=1.5cm and 2.5cm, on grid, auto]

% Styles
\tikzstyle{input} = [circle, draw=black, fill=tabgreen!30, thick, minimum size=8mm]
\tikzstyle{output} = [circle, draw=black, fill=tabolive!30, thick, minimum size=8mm]
%\tikzstyle{hidden} = [rectangle, text=black, text opacity=1, draw=black, fill=black!40, fill opacity=0.3, thick, minimum size=8mm, rounded corners=3pt]
\newcommand{\hidden}[4]{% #1=x, #2=y, #3=name, #4=label
  \node[inner sep=0pt, outer sep=0pt] (#3) at (#1,#2) {};
  % Draw the polygon
  \draw[fill=grey!20, draw=black]
    (#3)
    -- ++(0,0.2*1.5)
    -- ++(0.2*5,0)
    -- ++(0.2*2.5,-0.2*2)
    -- ++(-0.2*2.5,-0.2*2)
    -- ++(-0.2*5,0)
    -- cycle;
  % Define right tip coordinate manually
  \coordinate (#3-right) at ({#1 + 1.5}, {#2 - 0.1});
  % Label
  \node at ({#1 + 0.6}, {#2 - 0.1}) {\textbf{#4}};
}

% Nodes
\node[input] (I1) at (0,0) {$\hat{x}_1$};
\node[input] (I2) at (0,-1) {$\hat{x}_2$};

%\node[hidden] (H1) at (2.5,1) {$\mathrm{NNP}_1$};
%\node[hidden] (H2) at (2.5,0) {$\mathrm{NNP}_2$};
%\node[hidden] (H3) at (2.5,-1) {$\mathrm{NNP}_3$};
%\node[hidden] (H4) at (2.5,-2) {$\mathrm{NNP}_4$};

\hidden{1.75}{1.0}{H1}{$\mathrm{NNP}_1$};
\hidden{1.75}{0.0}{H2}{$\mathrm{NNP}_2$};
\hidden{1.75}{-1.0}{H3}{$\mathrm{NNP}_3$};
\hidden{1.75}{-2.0}{H4}{$\mathrm{NNP}_4$};

\node[output] (O1) at (5,0) {$\dot{\hat{x}}_1$};
\node[output] (O2) at (5,-1) {$\dot{\hat{x}}_2$};

% Connections from Input to Hidden
\foreach \i in {I1,I2} {
  \foreach \h in {H1,H2,H3,H4} {
    \draw[gray] (\i) -- (\h);
  }
}

% Connections from Hidden to Output
\draw[tabgreen!60] (H1-right) -- (O1);
\draw[tabblue!60] (H2-right) -- (O1);
\draw[tabgreen!60] (H3-right) -- (O2);
\draw[tabblue!60] (H4-right) -- (O2);

\node at (0,-3) (inputlabel) {input layer};
\node at (2.5,-3) (hiddenlabel) {NNP layer};
\node at (5,-3) (outputlabel) {output layer};

\node at (-0.5,1.5) {(b)};

\end{tikzpicture}

%% file: Monod_BiNODE.tikz
\begin{tikzpicture}[node distance=1.5cm and 2.5cm, on grid, auto]

% Styles
\tikzstyle{input} = [circle, draw=black, fill=tabgreen!30, thick, minimum size=8mm]
\tikzstyle{output} = [circle, draw=black, fill=tabolive!30, thick, minimum size=8mm]
%\tikzstyle{hidden} = [rectangle, text=black, text opacity=1, draw=black, fill=black!40, fill opacity=0.3, thick, minimum size=8mm, rounded corners=3pt]
\newcommand{\hidden}[4]{% #1=x, #2=y, #3=name, #4=label
  \node[inner sep=0pt, outer sep=0pt] (#3) at (#1,#2) {};
  % Draw the polygon
  \draw[fill=grey!20, draw=black]
    (#3)
    -- ++(0,0.2*1.5)
    -- ++(0.2*5,0)
    -- ++(0.2*2.5,-0.2*2)
    -- ++(-0.2*2.5,-0.2*2)
    -- ++(-0.2*5,0)
    -- cycle;
  % Define right tip coordinate manually
  \coordinate (#3-right) at ({#1 + 1.5}, {#2 - 0.1});
  % Label
  \node at ({#1 + 0.6}, {#2 - 0.1}) {\textbf{#4}};
}

% Nodes
\node[input] (I1) at (0,0.5) {$\hat{x}_1$};
\node[input] (I2) at (0,-0.5) {$\hat{x}_2$};

\hidden{1.75}{0.5}{H1}{$\mathrm{NNP}_1$};
\hidden{1.75}{-0.5}{H2}{$\mathrm{NNP}_2$};

\node[output] (O1) at (5,0.4) {$\dot{\hat{x}}_1$};
\node[output] (O2) at (5,-0.6) {$\dot{\hat{x}}_2$};

% Connections from Input to Hidden
\foreach \i in {I2} {
  \foreach \h in {H1} {
    \draw[gray] (\i) -- (\h);
  }
}
\foreach \i in {I1} {
  \foreach \h in {H1,H2} {
    \draw[gray] (\i) -- (\h);
  }
}

% Connections from Hidden to Output
\foreach \h in {H1,H2} {
  \foreach \o in {O1} {
    \draw[gray] (\h-right) -- (\o);
  }
}

\draw[gray] (H1-right) -- (O2);

\node at (0,-2.25) (inputlabel) {input layer};
\node at (2.5,-2.25) (hiddenlabel) {NNP layer};
\node at (5,-2.25) (outputlabel) {output layer};

%\node at (-0.5,1.5) {(a))};%

\end{tikzpicture}

%% file: LV_BiNODE.tikz
\begin{tikzpicture}[node distance=1.5cm and 2.5cm, on grid, auto]

% Styles
\tikzstyle{input} = [circle, draw=black, fill=tabgreen!30, thick, minimum size=8mm]
\tikzstyle{output} = [circle, draw=black, fill=tabolive!30, thick, minimum size=8mm]
%\tikzstyle{hidden} = [rectangle, text=black, text opacity=1, draw=black, fill=black!40, fill opacity=0.3, thick, minimum size=8mm, rounded corners=3pt]
\newcommand{\hidden}[4]{% #1=x, #2=y, #3=name, #4=label
  \node[inner sep=0pt, outer sep=0pt] (#3) at (#1,#2) {};
  % Draw the polygon
  \draw[fill=grey!20, draw=black]
    (#3)
    -- ++(0,0.2*1.5)
    -- ++(0.2*5,0)
    -- ++(0.2*2.5,-0.2*2)
    -- ++(-0.2*2.5,-0.2*2)
    -- ++(-0.2*5,0)
    -- cycle;
  % Define right tip coordinate manually
  \coordinate (#3-right) at ({#1 + 1.5}, {#2 - 0.1});
  % Label
  \node at ({#1 + 0.6}, {#2 - 0.1}) {\textbf{#4}};
}

% Nodes
\node[input] (I1) at (0,0.5) {$\hat{x}_1$};
\node[input] (I2) at (0,-0.5) {$\hat{x}_2$};

\hidden{1.75}{1.0}{H1}{$\mathrm{NNP}_1$};
\hidden{1.75}{0.0}{H2}{$\mathrm{NNP}_2$};
\hidden{1.75}{-1.0}{H3}{$\mathrm{NNP}_3$};
%\hidden{1.75}{-2.0}{H4}{$\mathrm{NNP}_4$};%

\node[output] (O1) at (5,0.5) {$\dot{\hat{x}}_1$};
\node[output] (O2) at (5,-0.5) {$\dot{\hat{x}}_2$};

% Connections from Input to Hidden
\foreach \i in {I1,I2} {
  \foreach \h in {H3} {
    \draw[gray] (\i) -- (\h);
  }
}

\draw[gray] (I1) -- (H1);
\draw[gray] (I2) -- (H2);

% Connections from Hidden to Output
\draw[gray] (H1-right) -- (O1);
\draw[gray] (H3-right) -- (O1);
\draw[gray] (H2-right) -- (O2);
\draw[gray] (H3-right) -- (O2);

\node at (0,-2.25) (inputlabel) {input layer};
\node at (2.5,-2.25) (hiddenlabel) {NNP layer};
\node at (5,-2.25) (outputlabel) {output layer};

%\node at (-0.5,1.5) {(a))};%

\end{tikzpicture}

%% file: refs.bib
@article{hossain2024biologically,
  title={Biologically informed {N}eural{ODE}s for genome-wide regulatory dynamics},
  author={Hossain, Intekhab and Fanfani, Viola and Fischer, Jonas and Quackenbush, John and Burkholz, Rebekka},
  journal={Genome Biology},
  volume={25},
  number={1},
  pages={127},
  year={2024},
  publisher={Springer}
}

@article{michaelis1913kinetik,
  title={{D}ie {K}inetik der {I}nvertinwirkung},
  author={Michaelis, Leonor and Menten, Maud L},
  journal={Biochemische Zeitschrift},
  volume={49},
  number={333-369},
  pages={352},
  year={1913},
  publisher={Berlin}
}

@book{henri1903lois,
  title={Lois g{\'e}n{\'e}rales de l'action des diastases},
  author={Henri, Victor},
  year={1903},
  publisher={Librairie Scientifique A. Hermann}
}

@article{henri2006theorie,
  title={Th{\'e}orie g{\'e}n{\'e}rale de l'action de quelques diastases par {V}ictor {H}enri [{CR} {A}cad. {S}ci. {P}aris 135 (1902) 916-919]},
  author={Henri, Victor},
  journal={Comptes Rendus Biologies},
  volume={329},
  number={1},
  pages={47--50},
  year={2006},
  publisher={Elsevier}
}

@article{cornish2015one,
  title={One hundred years of Michaelis--Menten kinetics},
  author={Cornish-Bowden, Athel},
  journal={Perspectives in Science},
  volume={4},
  pages={3--9},
  year={2015},
  publisher={Elsevier}
}

@article{waage1864studier,
  title={Studier over affiniteten},
  author={Waage, P and Guldberg, CM},
  journal={Forhandlinger i Videnskabs-selskabet i Christiania},
  volume={1},
  pages={35--45},
  year={1864}
}

@article{voit2015150,
  title={150 years of the mass action law},
  author={Voit, Eberhard O and Martens, Harald A and Omholt, Stig W},
  journal={PLOS Computational Biology},
  volume={11},
  number={1},
  pages={e1004012},
  year={2015},
  publisher={Public Library of Science San Francisco, USA}
}

@article{kermack1927contribution,
  title={A contribution to the mathematical theory of epidemics},
  author={Kermack, William Ogilvy and McKendrick, Anderson G},
  journal={Proceedings of the Royal Society of London. Series A, Containing Papers of a Mathematical and Physical Character},
  volume={115},
  number={772},
  pages={700--721},
  year={1927},
  publisher={The Royal Society London}
}

@book{volterra1927variazioni,
  title={Variazioni e fluttuazioni del numero d'individui in specie animali conviventi},
  author={Volterra, Vito},
  volume={2},
  year={1927},
  publisher={Societ{\'a} anonima tipografica" Leonardo da Vinci"}
}

@book{lotka1925elements,
  title={Elements of physical biology},
  author={Lotka, Alfred James},
  year={1925},
  publisher={Williams \& Wilkins}
}

@article{visser2003dynamic,
  title={Dynamic simulation and metabolic re-design of a branched pathway using linlog kinetics},
  author={Visser, Diana and Heijnen, Joseph J},
  journal={Metabolic engineering},
  volume={5},
  number={3},
  pages={164--176},
  year={2003},
  publisher={Elsevier}
}

@inproceedings{kacser1973control,
  title={The control of flux},
  author={Kacser, Henrik},
  booktitle={Symp Soc Exp Biol},
  volume={27},
  pages={65},
  year={1973}
}

@article{heinrich1974linear,
  title={A linear steady-state treatment of enzymatic chains: general properties, control and effector strength},
  author={Heinrich, Reinhart and Rapoport, Tom A},
  journal={European Journal of Biochemistry},
  volume={42},
  number={1},
  pages={89--95},
  year={1974},
  publisher={Wiley Online Library}
}

@article{savageau1969iochemical,
  title={Biochemical systems analysis: {I}. {S}ome mathematical properties of the rate law for the component enzymatic reactions},
  author={Savageau, Michael A},
  journal={Journal of Theoretical Biology},
  volume={25},
  number={3},
  pages={365--369},
  year={1969},
  publisher={Elsevier}
}

@article{savageau1969biochemical,
  title={Biochemical systems analysis: {II}. {T}he steady-state solutions for an n-pool system using a power-law approximation},
  author={Savageau, Michael A},
  journal={Journal of Theoretical Biology},
  volume={25},
  number={3},
  pages={370--379},
  year={1969},
  publisher={Elsevier}
}

@article{savageau1970biochemical,
  title={Biochemical systems analysis: {III}. {D}ynamic solutions using a power-law approximation},
  author={Savageau, Michael A},
  journal={Journal of Theoretical Biology},
  volume={26},
  number={2},
  pages={215--226},
  year={1970},
  publisher={Elsevier}
}

@article{voit2013biochemical,
  title={{B}iochemical {S}ystems {T}heory: {A} {R}eview},
  author={Voit, Eberhard O},
  journal={ISRN Biomathematics},
  volume={2013},
  year={2013},
  publisher={Hindawi}
}

@article{cleland1975partition,
  title={Partition analysis and concept of net rate constants as tools in enzyme kinetics},
  author={Cleland, WW},
  journal={Biochemistry},
  volume={14},
  number={14},
  pages={3220--3224},
  year={1975},
  publisher={ACS Publications}
}

@article{king1956schematic,
  title={A schematic method of deriving the rate laws for enzyme-catalyzed reactions},
  author={King, Edward L and Altman, Carl},
  journal={The Journal of Physical Chemistry},
  volume={60},
  number={10},
  pages={1375--1378},
  year={1956},
  publisher={ACS Publications}
}

@article{philipps2024universal,
  title={Universal differential equations for systems biology: Current state and open problems},
  author={Philipps, Maren and Schmid, Nina and Hasenauer, Jan},
  journal={bioRxiv},
  pages={2024--11},
  year={2024}
}

@article{de2025physiology,
  title={Physiology-informed regularisation enables training of universal differential equation systems for biological applications},
  author={de Rooij, Max and Erd{\H{o}}s, Bal{\'a}zs and van Riel, Natal AW and O’Donovan, Shauna D},
  journal={PLOS Computational Biology},
  volume={21},
  number={1},
  pages={e1012198},
  year={2025}
}

@article{ahmadi2023learning,
  title={Learning dynamical systems with side information},
  author={Ahmadi, Amir Ali and Khadir, Bachir El},
  journal={SIAM Review},
  volume={65},
  number={1},
  pages={183--223},
  year={2023}
}

@article{grigorian2024hybrid,
  title={A hybrid neural ordinary differential equation model of the cardiovascular system},
  author={Grigorian, Gevik and George, Sandip V and Lishak, Sam and Shipley, Rebecca J and Arridge, Simon},
  journal={Journal of the Royal Society Interface},
  volume={21},
  number={212},
  pages={20230710},
  year={2024}
}

@article{bottcher2025control,
  title={Control of medical digital twins with artificial neural networks},
  author={B{\"o}ttcher, Lucas and Fonseca, Luis L and Laubenbacher, Reinhard C},
  journal={Philosophical Transactions A},
  volume={383},
  number={2292},
  pages={20240228},
  year={2025}
}

@article{bottcher2026control,
  title={Control of dynamical systems with neural networks},
  author={B{\"o}ttcher, Lucas},
  journal={Nonlinear Dynamics},
  volume={114},
  number={2},
  pages={79},
  year={2026}
}

@article{hill1910possible,
  title={The possible effects of the aggregation of the molecules of hemoglobin on its dissociation curves},
  author={Hill, Archibald Vivian},
  journal={The Journal of Physiology},
  volume={40},
  pages={iv--vii},
  year={1910}
}

@article{gesztelyi2012hill,
  title={The {H}ill equation and the origin of quantitative pharmacology},
  author={Gesztelyi, Rudolf and Zsuga, Judit and Kemeny-Beke, Adam and Varga, Balazs and Juhasz, Bela and Tosaki, Arpad},
  journal={Archive for history of exact sciences},
  volume={66},
  pages={427--438},
  year={2012},
  publisher={Springer}
}

@inproceedings{DBLP:conf/l4dc/AhmadiK20,
  author       = {Amir Ali Ahmadi and
                  Bachir El Khadir},
  editor       = {Alexandre M. Bayen and
                  Ali Jadbabaie and
                  George J. Pappas and
                  Pablo A. Parrilo and
                  Benjamin Recht and
                  Claire J. Tomlin and
                  Melanie N. Zeilinger},
  title        = {Learning Dynamical Systems with Side Information},
  booktitle    = {Proceedings of the 2nd Annual Conference on Learning for Dynamics
                  and Control, {L4DC} 2020, Online Event, Berkeley, CA, USA, 11-12 June
                  2020},
  series       = {Proceedings of Machine Learning Research},
  volume       = {120},
  pages        = {718--727},
  publisher    = {{PMLR}},
  year         = {2020},
  url          = {http://proceedings.mlr.press/v120/ahmadi20a.html}
}

@article{briggs1925note,
  title={A note on the kinetics of enzyme action},
  author={Briggs, George Edward and Haldane, John Burdon Sanderson},
  journal={Biochemical Journal},
  volume={19},
  number={2},
  pages={338},
  year={1925}
}

@phdthesis{monod1958recherches,
  title={Recherches sur la croissance des cultures bact{\'e}riennes},
  author={Monod, Jacques},
  school={Universit{\'e} de Paris},
  year={1941}
}

@inproceedings{monod1950technique,
  title={La technique de culture continue, th $\{${\'e}$\}$ orie et applications},
  author={Monod, Jacques},
  booktitle={Annales de l'Institut Pasteur},
  volume={79},
  pages={390--410},
  year={1950}
}

@article{hofmeyr1997reversible,
  title={The reversible {H}ill equation: how to incorporate cooperative enzymes into metabolic models},
  author={Hofmeyr, Jan-Hendrik S and Cornish-Bowden, Hofmeyr},
  journal={Bioinformatics},
  volume={13},
  number={4},
  pages={377--385},
  year={1997},
  publisher={Oxford University Press}
}

@article{liebermeister2010modular,
  title={Modular rate laws for enzymatic reactions: thermodynamics, elasticities and implementation},
  author={Liebermeister, Wolfram and Uhlendorf, Jannis and Klipp, Edda},
  journal={Bioinformatics},
  volume={26},
  number={12},
  pages={1528--1534},
  year={2010},
  publisher={Oxford University Press}
}

@article{liebermeister2006bringing,
  title={Bringing metabolic networks to life: convenience rate law and thermodynamic constraints},
  author={Liebermeister, Wolfram and Klipp, Edda},
  journal={Theoretical Biology and Medical Modelling},
  volume={3},
  pages={1--13},
  year={2006},
  publisher={Springer}
}

@article{sorribas2007cooperativity,
  title={Cooperativity and saturation in biochemical networks: a saturable formalism using {T}aylor series approximations},
  author={Sorribas, Albert and Hern{\'a}ndez-Bermejo, Benito and Vilaprinyo, Ester and Alves, Rui},
  journal={Biotechnology and Bioengineering},
  volume={97},
  number={5},
  pages={1259--1277},
  year={2007},
  publisher={Wiley Online Library}
}

@article{cs1959components,
  title={The components of predation as revealed by a study of small mammal predation of the {E}uropean pine sawfly},
  author={C S Holling},
  journal={The Canadian Entomologist},
  volume={91},
  pages={293--320},
  year={1959}
}

@article{holling1959some,
  title={Some characteristics of simple types of predation and parasitism},
  author={Holling, Crawford S},
  journal={The Canadian Entomologist},
  volume={91},
  number={7},
  pages={385--398},
  year={1959},
  publisher={Cambridge University Press}
}

@article{dawes2013derivation,
  title={A derivation of {H}olling's type {I}, {II} and {III} functional responses in predator--prey systems},
  author={Dawes, JHP and Souza, MO3046076},
  journal={Journal of Theoretical Biology},
  volume={327},
  pages={11--22},
  year={2013},
  publisher={Elsevier}
}

@article{wylie2021uniformly,
  title={Uniformly accurate nonlinear transmission rate models arising from disease spread through pair contacts},
  author={Wylie, Jonathan and Chou, Tom},
  journal={Physical Review E},
  volume={103},
  number={3},
  pages={032306},
  year={2021},
  publisher={APS}
}

@article{dietrich2023learning,
  title={Learning effective stochastic differential equations from microscopic simulations: Linking stochastic numerics to deep learning},
  author={Dietrich, Felix and Makeev, Alexei and Kevrekidis, George and Evangelou, Nikolaos and Bertalan, Tom and Reich, Sebastian and Kevrekidis, Ioannis G},
  journal={Chaos: An Interdisciplinary Journal of Nonlinear Science},
  volume={33},
  number={2},
  year={2023},
  publisher={AIP Publishing}
}

@article{zhang2025reconstructing,
  title={Reconstructing noisy gene regulation dynamics using extrinsic-noise-driven neural stochastic differential equations},
  author={Zhang, Jiancheng and Li, Xiangting and Guo, Xiaolu and You, Zhaoyi and B{\"o}ttcher, Lucas and Mogilner, Alex and Hoffmann, Alexander and Chou, Tom and Xia, Mingtao},
  journal={PLOS Computational Biology},
  volume={21},
  number={9},
  pages={e1013462},
  year={2025}
 }

@article{fonseca2025optimal,
  title={Optimal control of agent-based models via surrogate modeling},
  author={Fonseca, Luis L and B{\"o}ttcher, Lucas and Mehrad, Borna and Laubenbacher, Reinhard C},
  journal={PLOS Computational Biology},
  volume={21},
  number={1},
  pages={e1012138},
  year={2025},
  publisher={Public Library of Science San Francisco, CA USA}
}

@article{fronk2023interpretable,
  title={Interpretable polynomial neural ordinary differential equations},
  author={Fronk, Colby and Petzold, Linda},
  journal={Chaos: An Interdisciplinary Journal of Nonlinear Science},
  volume={33},
  number={4},
  year={2023},
  publisher={AIP Publishing}
}

@article{thöni2025modeling,
  title={Modeling Chemical Reaction Networks Using Neural Ordinary Differential Equations},
  author={Th\"oni, Anna CM and Robinson, William E and Bachrach, Yoram and Huck, Wilhelm TS and Kachman, Tal},
  journal={Journal of Chemical Information and Modeling},
  year={2025},
  publisher={ACS Publications}
}

@article{hornik1989multilayer,
  title={Multilayer feedforward networks are universal approximators},
  author={Hornik, Kurt and Stinchcombe, Maxwell and White, Halbert},
  journal={Neural networks},
  volume={2},
  number={5},
  pages={359--366},
  year={1989},
  publisher={Elsevier}
}

@inproceedings{DBLP:conf/nips/LuPWH017,
  author       = {Zhou Lu and
                  Hongming Pu and
                  Feicheng Wang and
                  Zhiqiang Hu and
                  Liwei Wang},
  editor       = {Isabelle Guyon and
                  Ulrike von Luxburg and
                  Samy Bengio and
                  Hanna M. Wallach and
                  Rob Fergus and
                  S. V. N. Vishwanathan and
                  Roman Garnett},
  title        = {The Expressive Power of Neural Networks: {A} View from the Width},
  booktitle    = {Advances in Neural Information Processing Systems 30: Annual Conference
                  on Neural Information Processing Systems 2017, December 4-9, 2017,
                  Long Beach, CA, {USA}},
  pages        = {6231--6239},
  year         = {2017},
  url          = {https://proceedings.neurips.cc/paper/2017/hash/32cbf687880eb1674a07bf717761dd3a-Abstract.html}
}

@inproceedings{DBLP:conf/colt/Telgarsky16,
  author       = {Matus Telgarsky},
  editor       = {Vitaly Feldman and
                  Alexander Rakhlin and
                  Ohad Shamir},
  title        = {Benefits of depth in neural networks},
  booktitle    = {Proceedings of the 29th Conference on Learning Theory, {COLT} 2016,
                  New York, USA, June 23-26, 2016},
  series       = {{JMLR} Workshop and Conference Proceedings},
  volume       = {49},
  pages        = {1517--1539},
  publisher    = {JMLR.org},
  year         = {2016},
  url          = {http://proceedings.mlr.press/v49/telgarsky16.html}
}

@article{cornish2013origins,
  title={The origins of enzyme kinetics},
  author={Cornish-Bowden, Athel},
  journal={FEBS letters},
  volume={587},
  number={17},
  pages={2725--2730},
  year={2013},
  publisher={Elsevier}
}

@article{haldane1930enzymes,
  title={Enzymes longmans},
  author={Haldane, JBS},
  journal={Green and Co, UK},
  volume={7},
  year={1930}
}

@inproceedings{DBLP:conf/acssc/DouglasY18,
  author       = {Scott C. Douglas and
                  Jiutian Yu},
  editor       = {Michael B. Matthews},
  title        = {Why {RELU} Units Sometimes Die: Analysis of Single-Unit Error Backpropagation
                  in Neural Networks},
  booktitle    = {52nd Asilomar Conference on Signals, Systems, and Computers, {ACSSC}
                  2018, Pacific Grove, CA, USA, October 28-31, 2018},
  pages        = {864--868},
  publisher    = {{IEEE}},
  year         = {2018},
  url          = {https://doi.org/10.1109/ACSSC.2018.8645556},
  doi          = {10.1109/ACSSC.2018.8645556},
}

@article{acon2021myc,
  title={MYC dosage compensation is mediated by miRNA-transcription factor interactions in aneuploid cancer},
  author={Ac{\'o}n, ManSai and Gei{\ss}, Carsten and Torres-Calvo, Jorge and Bravo-Estupi{\~n}an, Diana and Oviedo, Guillermo and Arias-Arias, Jorge L and Rojas-Matey, Luis A and Edwin, Baez and V{\'a}squez-Vargas, Gloriana and Oses-Vargas, Yendry and others},
  journal={IScience},
  volume={24},
  number={12},
  year={2021},
  publisher={Elsevier}
}

@article{andrews1968mathematical,
  title={A mathematical model for the continuous culture of microorganisms utilizing inhibitory substrates},
  author={Andrews, John F},
  journal={Biotechnology and Bioengineering},
  volume={10},
  number={6},
  pages={707--723},
  year={1968},
  publisher={Wiley Online Library}
}

@article{muloiwa2020comparison,
  title={Comparison of unstructured kinetic bacterial growth models},
  author={Muloiwa, Mpho and Nyende-Byakika, Stephen and Dinka, Megersa},
  journal={South African Journal of Chemical Engineering},
  volume={33},
  pages={141--150},
  year={2020},
  publisher={Elsevier}
}

@book{haldane1965enzymes,
  title={Enzymes},
  author={Haldane, JBS},
  year={1965},
  publisher={MIT Press},
  address={Cambridge, MA, USA}
}

@book{moser1958dynamics,
  title={The dynamics of bacterial populations maintained in the chemostat.},
  author={Moser, Hermann},
  year={1958},
  publisher={Carnegie Institution of Washington}
}

@article{AI-Aristotle,
    doi = {10.1371/journal.pcbi.1011916},
    author = {Ahmadi Daryakenari, Nazanin AND De Florio, Mario AND Shukla, Khemraj AND Karniadakis, George Em},
    journal = {PLOS Computational Biology},
    publisher = {Public Library of Science},
    title = {{AI}-{A}ristotle: A physics-informed framework for systems biology gray-box identification},
    year = {2024},
    month = {03},
    volume = {20},
    url = {https://doi.org/10.1371/journal.pcbi.1011916},
    pages = {1-33},
    number = {3},

}

@article{schiassi2021X-TFC,
  title={Extreme theory of functional connections: A fast physics-informed neural network method for solving ordinary and partial differential equations},
  author={Schiassi, Enrico and Furfaro, Roberto and Leake, Carl and De Florio, Mario and Johnston, Hunter and Mortari, Daniele},
  journal={Neurocomputing},
  volume={457},
  pages={334--356},
  year={2021},
  publisher={Elsevier}
}

@article{de2022PINNs,
  title={Physics-informed neural networks and functional interpolation for stiff chemical kinetics},
  author={De Florio, Mario and Schiassi, Enrico and Furfaro, Roberto},
  journal={Chaos: An Interdisciplinary Journal of Nonlinear Science},
  volume={32},
  number={6},
  year={2022},
  publisher={AIP Publishing}
}

@article{virgolin2022symbolic,
  title={Symbolic regression is {NP}-hard},
  author={Virgolin, Marco and Pissis, Solon P},
  journal={arXiv preprint arXiv:2207.01018},
  year={2022}
}

@misc{gplearn,
  author = {Trevor Stephens},
  title = {{gplearn}: Genetic Programming in Python with a scikit-learn inspired and compatible API},
  year = {2026},
  url = {https://github.com/trevorstephens/gplearn},
  note = {Version 0.4.3},
}

@book{barnes2011mathematical,
  title={Mathematical modelling with case studies: a differential equations approach using Maple and {MATLAB}},
  author={Barnes, Belinda and Fulford, G R},
  year={2011},
  publisher={Chapman and Hall/CRC}
}

@inproceedings{DBLP:conf/icml/AmosXK17,
  author       = {Brandon Amos and
                  Lei Xu and
                  J. Zico Kolter},
  editor       = {Doina Precup and
                  Yee Whye Teh},
  title        = {Input Convex Neural Networks},
  booktitle    = {Proceedings of the 34th International Conference on Machine Learning,
                  {ICML} 2017, Sydney, NSW, Australia, 6-11 August 2017},
  series       = {Proceedings of Machine Learning Research},
  volume       = {70},
  pages        = {146--155},
  publisher    = {{PMLR}},
  year         = {2017},
  url          = {http://proceedings.mlr.press/v70/amos17b.html}
}

@article{reardon2000biodegradation,
  title={Biodegradation kinetics of benzene, toluene, and phenol as single and mixed substrates for {P}seudomonas putida {F1}},
  author={Reardon, Kenneth F and Mosteller, Douglas C and Bull Rogers, Julia D},
  journal={Biotechnology and Bioengineering},
  volume={69},
  number={4},
  pages={385--400},
  year={2000}
}

@article{sturis1991computer,
  title={Computer model for mechanisms underlying ultradian oscillations of insulin and glucose},
  author={Sturis, Jeppe and Polonsky, Kenneth S and Mosekilde, Erik and Van Cauter, Eve},
  journal={American Journal of Physiology-Endocrinology and Metabolism},
  volume={260},
  number={5},
  pages={E801--E809},
  year={1991}
}

@article{desilva2020,
doi = {10.21105/joss.02104},
url = {https://doi.org/10.21105/joss.02104},
year = {2020},
publisher = {The Open Journal},
volume = {5},
number = {49},
pages = {2104},
author = {Brian de Silva and Kathleen Champion and Markus Quade and Jean-Christophe Loiseau and J. Kutz and Steven Brunton},
title = {{PySINDy}: A {P}ython package for the sparse identification of nonlinear dynamical systems from data},
journal = {Journal of Open Source Software}
}

@article{Kaptanoglu2022,
doi = {10.21105/joss.03994},
url = {https://doi.org/10.21105/joss.03994},
year = {2022},
publisher = {The Open Journal},
volume = {7},
number = {69},
pages = {3994},
author = {Alan A. Kaptanoglu and Brian M. de Silva and Urban Fasel and Kadierdan Kaheman and Andy J. Goldschmidt and Jared Callaham and Charles B. Delahunt and Zachary G. Nicolaou and Kathleen Champion and Jean-Christophe Loiseau and J. Nathan Kutz and Steven L. Brunton},
title = {{PySINDy}: A comprehensive Python package for robust sparse system identification},
journal = {Journal of Open Source Software}
}

@article{SINDy2016,
author = {Steven L. Brunton  and Joshua L. Proctor  and J. Nathan Kutz },
title = {Discovering governing equations from data by sparse identification of nonlinear dynamical systems},
journal = {Proceedings of the National Academy of Sciences},
volume = {113},
number = {15},
pages = {3932-3937},
year = {2016},
doi = {10.1073/pnas.1517384113},
}

@software{Kaptanoglu_pysindy,
author = {Kaptanoglu, Alan and Stevens-Haas, Jacob and Champion, Kathleen and de Silva, Brian and Quade, Markus},
license = {MIT},
title = {{PySINDy}},
url = {https://github.com/dynamicslab/pysindy}
}

@article{fonseca2011complex,
  title={Complex coordination of multi-scale cellular responses to environmental stress},
  author={Fonseca, Lu{\'\i}s L and Sanchez, Claudia and Santos, Helena and Voit, Eberhard O},
  journal={Molecular BioSystems},
  volume={7},
  number={3},
  pages={731--741},
  year={2011},
  publisher={The Royal Society of Chemistry Oxford, UK}
}
